\documentclass[a4paper,twoside,12pt]{amsart}

\usepackage{amsmath, amsfonts, amssymb}
\usepackage{graphicx}
\usepackage{amsthm}
\usepackage{amsmath}
\usepackage[ansinew]{inputenc}
\usepackage{xy}

\usepackage{amssymb}
\usepackage{amsmath}
\usepackage{epsfig}

\newtheorem{theorem}{Theorem}[section]
\newtheorem{definition}[theorem]{Definition}
\newtheorem{lemma}[theorem]{Lemma}
\newtheorem{proposition}[theorem]{Proposition}
\newtheorem{corollary}[theorem]{Corollary}

\newtheorem{remark}[theorem]{Remark}

\newtheorem{example}[theorem]{Example}

\def\a{\alpha}

\def\b{\beta}
\def\C{\mathbf C}
\def\D{\Delta}
\def\d{\delta}

\def\f{\phi}
\def\g{\gamma}

\def\l{\lambda}

\def\O{\mbox{\rm orb}}

\def\p{\pi}

\def\P{{\mathbf P}}

\def\Q{\mathbf Q}
\def\R{\mathbf R}
\def\r{\rho}

\def\s{\sigma}
\def\t{\tau}

\def\w{\omega}
\def\Z{\mathbf Z}
\def\z{\zeta}

\title{Quasi ordinary singularities, essential divisors and Poincaré series }
\author{P.D. Gonz\'alez P\'erez}

\address{Departamento de Algebra. Facultad de Ciencias Matem\'aticas. Universidad Complutense de Madrid.
Plaza de las Ciencias 3. 28040. Madrid. Spain.}

\email{pgonzalez@mat.ucm.es}
\thanks{Gonz\'alez P\'erez is supported by {\em Programa Ram\'on y Cajal} and  MTM2004-08080-C02-01 grants of {\em Ministerio de Educaci\'on y Ciencia}, Spain.}

\author{F. Hernando}
\address{Departamento de Algebra, Geometría y Topología,
Universidad de Valladolid. Valladolid. Spain.}

\email{hernando@agt.uva.es}
\thanks{Hernando is supported by grants  MTM2004-010958 of {\em Ministerio de Educaci\'on y
Ciencia} and VA068/04 of {\em Junta de Castilla y León}, Spain.}


\begin{document}

\begin{abstract}
We define Poincaré series associated to a toric or analytically
irreducible quasi-ordinary hypersurface singularity, $(S,0)$, by a
finite sequence of monomial valuations, such that at least one of
them is centered at the origin $0$. This involves the definition
of a multigraded ring associated to the analytic algebra of the
singularity by the sequence of valuations. We prove that the
Poincaré series is a rational function with integer coefficients,
which  can be defined also as an integral with respect of the
Euler characteristic, over the projectivization of the analytic
algebra of the singularity, of a function defined by the
valuations. In particular, the Poincaré series associated to the
set of divisorial valuations corresponding to the essential
divisors, considered both over the singular locus and over the
point $0$, is an analytic invariant of the singularity. In the
quasi-ordinary hypersurface case we prove that this Poincaré
series determines and it is determined by the normalized sequence
of characteristic monomials. These monomials in the analytic case
define a complete invariant of the embedded topological type of
the hypersurface singularity.
\end{abstract}

 \maketitle


\section{Introduction}
\subsection{Poincaré series of singularities}

In recent times there has been many papers which analyze the
structure of a germ of complex analytic (or algebroid) singularity
$(S,0)$ of pure dimension $d$ by attaching to it a suitable notion
of Poincaré series (see for instance \cite{GDC},
\cite{C-D-G-curvas},\cite{C-D-G-superficies-racionales},
\cite{C-H-R-surface-singularities}, \cite{Ann-toricas},
\cite{E-S-quasihomogeneous} and
\cite{E-S-quasihomogeneous-complete-intersection}).

The basic case, when $(S,0)$ is an irreducible complex analytic
plane curve germ in $(\C^2, 0)$ defined by an equation $f(x, y)
=0$, was analyzed by Gusein-Zade, Delgado and Campillo in
\cite{GDC}. They show by a direct computation that if $\Gamma
\subset \Z_{\geq 0}$ is the semigroup of the curve, the Poincaré
series $\sum_{\g \in \Gamma} t^{\g}$ coincides with the zeta
function of geometric monodromy of the function $f: (\C^2, 0)
\rightarrow (\C,0)$. Relevant to this case is the idea of studying
a germ of singularity by analyzing the properties of the
filtrations of the analytic algebra $\mathcal{O}_S$ and the
associated graded rings, in particular using valuations (see for
instance the work of Lejeune and Teissier \cite{LejT}). Lejeune's
observation, applied to this basic case, was that the semigroup
$\Gamma$ is determined by the graded ring associated to
$\mathcal{O}_S$, by the filtration induced by the maximal ideal of
its integral closure (see \cite{Teissier} page 161).

Given a finite sequence $\underline{\nu} = (\nu_1,\ldots,\nu_r)$
of divisorial valuations of the analytic algebra $\mathcal{O}_S$,
Campillo, Delgado and Gusein-Zade defined in several context a
multi-index filtration of the analytic algebra $\mathcal{O}_S$,
associated to a finite sequence $\underline{\nu} =
(\nu_1,\ldots,\nu_r)$ of divisorial valuations. An ideal of the
filtration, which have index $\underline{n} =(n_1, \dots, n_r) \in
\Z^r$, is of the form $I_{\underline{n}}=\{g\in \mathcal{O}_S \mid
\underline{\nu} (g)  \geq \underline{n} \}$, where $\geq$ means
coordinate-wise and $\underline{\nu} (g) = (\nu_1 (g) , \dots,
\nu_r(g))$. If all the valuations in the sequence are centered at
the origin $0 \in S$ then
 the dimension
 \begin{equation} \label{cdim}
c_{\underline{n}}:=dim_{\mathbb{C}}(I_{\underline{n}}/I_{\underline{n}+\underline{1}})
\end{equation}
are finite for all $ \underline{n} \in \Z^{r}$. These authors
define then the  associated Poicaré series as the formal power
series in the indeterminates $\underline{t} = (t_1, \dots, t_r)$
given by:
\begin{equation}\label{def-C-D-G}
P =(\sum_{\underline{n}\in\mathbb{Z}^r} c_{\underline{n}}
\underline{t}^{\underline{n}})\frac{\prod_{i=1}^r
(t_i-1)}{\underline{t}^{\underline{1}}-1} \, \, \mbox{ where }
\underline{t}^{\underline{n}} = t_1^{n_1} \dots t_r^{n_r} \mbox{
and }   {\underline{1}} = (1,\dots, 1) .
\end{equation}
For technical reasons this definition is more convenient than
$
P' =\sum_{\underline{n}\in\mathbb{Z}_{\geq 0}^r} c_{\underline{n}}
\underline{t}^{\underline{n}}.
$

They actually compute these series, for instance, when $(S,0)$ is
a germ of plane curve using the valuations  defined by the
parametrizations of its irreducible components, or when $(S,0)$ is
a germ of rational surface singularity using the divisorial
valuations of the irreducible components of the exceptional
divisor of its minimal embedded resolution. The computation passes
through a suitable notion of integration with respect to the Euler
characteristic, over the projectivization of the algebra of the
singularity, of a function defined by the valuations.
 The result in the case of
plane curves coincides with the Alexander polynomial, which is a
topological invariant. In  the case of rational surface
singularities the Poincaré series is a rational function in
$\underline{t}$, with denominator determined explicitly by the
intersection matrix of the exceptional divisor in the minimal
resolution of the surface singularity, which is also a topological
invariant.

More generally if  $(S,0)$ is a germ of normal surface singularity
Cutkosky, Herzog and Reguera associated to a fixed resolution of
singularities of $(S,0)$ a Poincaré series defined by
\begin{equation}\label{def-usando-longitudes}
Q=\sum_{\underline{n}\in\mathbb{Z}_{\geq 0}^r} l_{\underline{n}}
\underline{t}^{\underline{n}}, \mbox{ where } \quad
l_{\underline{n}}= \mbox{length}( \mathcal{O}_S
/I_{\underline{n}}),
\end{equation}
where $I_{\underline{n}}$ is defined as before with respect to the
divisorial valuations associated to the irreducible components of
the exceptional divisor. This series, which is studied in a
general algebraic setting, even if the residue field of the local
algebra $\mathcal{O}_S$ is of positive characteristic, is not in
general a rational function of $\underline{t}$. In the particular
case of rational surface singularities this series is rational and
it is related with the Poincaré series (\ref{def-C-D-G}) (see
\cite{C-H-R-surface-singularities} and
\cite{C-D-G-superficies-racionales}). Recent results 
on Poincaré series of surface singularities are announced 
in \cite{Nemethi-series}.

Lemahieu studies the Poincaré series associated to a toric
singularity by a sequence of monomial valuations centered at the
origin using definition (\ref{def-C-D-G}) and finds out a formula
in terms of integration with respect to Euler characteristic over
the projectivization of the ring of functions \cite{Ann-toricas}.

\subsection{Essential divisors and the Nash map}

The divisors corresponding to the irreducible components of the
minimal embedded resolution of a germ of surface singularity
provide the basic example of {\em essential divisors}.  Roughly
speaking, essential is a property of the associated divisorial
valuation which must hold in every resolution of singularities.
This notion is much more interesting if the dimension of the
variety is $\geq 3$, since we do not dispose then of a unique
minimal resolution of singularities. We say that the divisorial
valuation associated to a essential divisor is also {\em
essential}.

Essential divisors appear naturally in relation with the Nash map.
Nash introduced a map from the set of families of arcs passing
through the singular locus of a variety to the set of essential
divisors over the singular locus. This map is called the Nash map.
We can also consider a similarly defined map from the set of
families of arcs passing through a fixed singular point and the
divisors over this point which are essential. This map is called
the local Nash map. If we have a variety with an isolated singular
locus both maps coincide. Both maps are injective (see
\cite{Nash}, \cite{IK}). Nash's question was on the surjectivity
of these maps. Ishii and Kollár showed a four dimensional isolated
hypersurface singularity for which the Nash map is not surjective
(see \cite{IK}). Nash problem is still open in dimensions two and
three and it is of particular interest to understand which are the
singularities such that the associated Nash map is bijective.

With this view-point it is very natural to consider an invariant
of the germ $(S,0)$ defined by the Poincaré series associated to
$(S,0)$ by the finite sequence of essential divisorial valuations.
If $(S,0)$ is an isolated singularity then we can define the
Poincaré series using the multi-index filtration and definition
(\ref{def-C-D-G}). Notice that  if $(S,0)$ is a non isolated
singularity it is possible that no essential valuation over the
singular locus of $S$ is centered over $0$, for instance if $d=2$
and $S$ is the total space of an equisingular family of plane
branches (the preimage of the singular locus in the normalization,
which in this case is a smooth surface, is the only essential
divisor). However, if one of the valuations is not centered at $0$
we cannot use definition (\ref{def-C-D-G}) to define the Poincaré
series since the dimensions (\ref{cdim}) become infinite. On the
other hand, considering only essential divisors over a fixed
singular point of $S$, does not take in to account the variety of
information hidden in the singular locus, when the singularity $S$
is not isolated.

We develop this program, extending to this case the definition of
the Poincaré series, when $(S,0)$ is an toric singularity or an
irreducible germ of quasi-ordinary hypersurface singularity, two
classes or singularities which are not isolated in general.

\subsection{Quasi-ordinary singularities and equisingularity} \label{equisingular}

Quasi-ordinary singularities arise classically in Jung's approach
to analyze surface singularities by using a finite projection to a
smooth surface (see \cite{Jung}). We say that a germ  $(S,0)
\subset (\C^{d+1}, 0) $ of complex analytic hypersurface with
equation $f=0$  is {\it quasi-ordinary} if there exists a finite
morphism $ \p: (S, 0) \rightarrow (\C^d,0)$ (called a {\it
quasi-ordinary projection\index{quasi-ordinary projection}}) such
that the discriminant locus is a germ of normal crossing divisor.
We will suppose from now on that the germ $(S,0)$ is analytically
irreducible, for simplicity. Quasi-ordinary hypersurface
singularities are parameterized by {\it quasi-ordinary branches},
certain class of fractional power series in several variables
having a finite set of {\it distinguished} or {\it characteristic
monomials}, which generalize
  the {\it characteristic pairs} associated to a plane branch
(see \cite{Zariski}). Since different quasi-ordinary projections
may define different sequence of characteristic monomials Lipman
showed that any quasi-ordinary hypersurface admits a {\em
normalized} quasi-ordinary projection, a technical condition which
adapts to this case the notion of transversal projection. Then he
proved in the embedded surface case, using a detailed analysis of
monoidal and quadratic transforms, that the sequence of normalized
characteristic monomials, which in principle depends  on the
corresponding quasi-ordinary projection, is in fact an analytic
invariant of the singularity
 (see \cite{Lipman65}, see also \cite{Luengo}).
  Lipman and Gau studied the
topological invariants of a germ of quasi-ordinary singularity.
Recall that two hypersurface germs, $(S_i, 0) \subset (\C^{d+1},
0)$ for $i=1,2$,  have the same embedded topological type if there
exists a germ homeomorphism $(\C^{d+1}, S_1, 0 ) \rightarrow
(\C^{d+1}, S_{2}, 0)$  of triples. Two irreducible germs of
$d$-dimensional quasi-ordinary hypersurface $(S_i, 0) \subset
(\C^{d+1},0)$ for $i=1,2$, have the same embedded topological type
if and only if both have the same sequence of normalized
characteristic monomials (see \cite{Lipman88}, \cite{Gau} and
\cite{Gau-T}).

Suppose that the reduced discriminant locus of the quasi-ordinary
projection $\p$ has equation $x_1 \cdots x_{c} =0$. Then we can
define for $I \subset \{ 1, \dots, c\}$ the sets: $ W_I :=  \{
(x_1, \dots, x_d) \mid x_i = 0, \, i\in I \}$, $W_I^\circ:= W_I -
\bigcup _{I\subset J} W_J$, $Z_I := \p^{-1} (W_I) \subset S$ and
$Z_I^\circ:= \p^{-1} ( W_I^\circ)$. Lipman showed that the sets
$\{ Z_I^\circ \}_I$ define an equimultiple complex analytic
stratification of $(S,0)$. As a corollary he proved that the
characteristic monomials determine the singular locus of $(S,0)$
(see  Theorem 7.3, \cite{Lipman88} and Theorem \ref{singular}). We
have that the restriction of $\p$ to $W = Z_{1, \dots, c} $
defines an isomorphism with $\{0 \} \times \C^{d-c}$. We  can
consider then $S$ as the total space of a family, $S_t$ for $t \in
(\C^{d-c},0)$, of hypersurfaces of dimension $c$, which are given
with a quasi-ordinary projection $\p_t : S_t \rightarrow \C^c$ for
each fiber. Gau and Lipman results imply that $(S,0)$ is
topologically equisingular along $W$, i.e., at each point of $W$
all fibers $S_t$ have the same topological type. Ban showed indeed
that $\{ Z_I^\circ \}_I$ defines a {\em Whitney stratification} of
$(S,0)$ and hence $(S,0)$ is {\em differentially equisingular}
along $(W,0)$, i.e., the pair $(S, W)$ satisfies the Whitney
conditions at $0$; conversely, if $(S,0)$ is differentially
equisingular along a smooth subgerm $(W,0)$  and the fibers $S_t$
are quasi-ordinary hypersurfaces of dimension $c$ for the same
projection, then all the fibers $S_t$
 have the same
embedded topological type (see \cite{Ban}). The integer $c$ is
called in \cite{patrick} the {\em equisingular dimension of
$(S,0)$}, since $(S,0)$ may be viewed as an equisingular
deformation of a quasi-ordinary hypersurface germ of dimension
$c$, but not smaller. See \cite{Ban} and also Lipman's work
\cite{Lipman97}.

Quasi-ordinary singularities play also an important role in the
comprehension of the different notions of  {\em equisingularity},
in particular with respect to {\em Zariski equisingularity} and
{\em equiresolution} (see \cite{Lipman97}). Villamayor proved that
if $\mathcal{W}$ is a nonsingular irreducible subvariety in an
algebraic hypersurface $\mathcal{S}$ inside a smooth algebraic
variety, then  if $\mathcal{S}$ is equisingular along
$\mathcal{W}$ in the sense of Zariski, then $\mathcal{S}$ is {\em
weakly equiresolvable} along $\mathcal{W}$ (see \cite{Vi}). His
proof passes through a construction of embedded resolution of
quasi-ordinary hypersurfaces, which is not constructed by blowing
up smooth centers. Embedded resolutions of a quasi-ordinary
hypersurface $(S,0)$ are also constructed,  using toric
modifications, from the associated characteristic monomials in
\cite{Pedro-resolucion-torica}. In the irreducible case this
method provides, after a suitable re-embedding of the germ in a
larger affine space, an embedded resolution with only one toric
modification, which desingularizes simultaneously the
quasi-ordinary germ $(S,0)$ and an affine toric variety $Z$, which
can be seen respectively as a generic fiber and special fiber of a
faithfully flat {\em equisingular} deformation. This approach
generalizes results of Teissier and Goldin for plane branches in
\cite{Rebeca}. The toric variety $Z $ is equal to $\mbox{Spec}
\C[\Gamma_\p]$ where the semigroup $\Gamma_\p$, defined in terms
of the characteristic monomials, was associated to the
quasi-ordinary projection $\p$ by Kiyek and Micus in \cite{Micus},
\cite{kiyek} and also in \cite{Pedro-Tesis},
\cite{Pedro-semigrupo} independently (the reference \cite{kiyek}
was not known by the first author at the time of writing these
works).

\subsection{The notion of semigroup of a quasi-ordinary
hypersurface}

Several papers have been made in the effort to settle down an
intrinsic notion of semigroup of an irreducible quasi-ordinary
hypersurface, which generalize to some extent the properties of
the semigroups of plane branches.   These papers and results are a
partial motivation for this paper.

- Kiyek and Micus in  \cite{kiyek} proved in the embedded surface
case that the semigroup $\Gamma_\p$ determines the  normalized
characteristic monomials, giving another proof of Lipman's result.

- In \cite{Pedro-Tesis} it was shown that the semigroups
$\Gamma_\p$ and $\Gamma_{\p'}$ associated respectively to a
quasi-ordinary projection $\p$ and normalized quasi-ordinary
projection $\p'$, defined from $\p$ by Lipman, are isomorphic. It
follows from the characterization of Lipman and Gau that the
isomorphy type of this semigroup is a complete invariant on the
embedded topological type (see \cite{Pedro-semigrupo}).

- In the embedded surface case Popescu-Pampu proved that the
semigroup $\Gamma_\p$ is an analytic invariant of the germ $f=0$,
by using the structure of the minimal embedded resolution of the
surface (see \cite{PPP02}).

- Popescu-Pampu  showed recently the analytical invariance of the
sequence of normalized characteristic monomials associated to
quasi-ordinary projection of $(S,0) \subset (\C^{d+1},0)$,
 without passing through
their embedded topological invariance. The argument involves the
definition of a canonical semigroup $\Gamma' (S)$ associated to
the germ $(S,0)$ as follows. First, he constructed a canonical
morphism of germs
\[
\theta : (\tilde{S}, P) \rightarrow (S,0),
\]
with source a smooth germ of dimension $d$, together with a
canonical germ of normal crossing divisor $\mathcal{H}$ at the
point $P$. The morphism $\theta$ and the divisor $\mathcal{H}$ are
constructed by using the structure of the pair $(S, \mbox{\rm
Sing} S)$ as follows: in order to obtain a smooth source the
normalization map $n: \bar{S} \rightarrow S$ is composed with a
finite canonical morphism, $\mu: S' \rightarrow \bar{S}$, defined
intrinsically in topological terms and called the {\em orbifold
mapping}. Then $(S',0) \cong (\C^d,0)$ and choosing suitable local
coordinates we have that $(n \circ \mu)^{-1} (\mbox{\rm Sing} S)$
is a union of finite number of coordinate subspaces. The map
$\eta: \tilde{S} \rightarrow S'$ is a sequence of blow ups with
smooth centers determined in combinatorial terms by $(n \circ
\mu)^{-1} (\mbox{\rm Sing} S)$. Then $\theta:= n \circ \mu \circ
\eta$ and this data defines a germ of normal crossing divisor germ
$\mathcal{H} := \theta^{-1}(\mbox{\rm Sing} S)$ at a point $P \in
\tilde{S}$ with $c'$ irreducible components, which correspond to
the first $c'$ elements  of a system   $u_1, \dots, u_d$ of local
coordinates at $P$. The subset of functions
\[
\mathcal{D} := \{ h \in \mathcal{O}_S \mid \theta^* (h) =
u_1^{m_1(h)} \cdots u_{c'} ^{m_{c'} (h)} \epsilon_h \mbox{ where }
\epsilon_h \in \C[[u_1, \dots, u_d]] \mbox{ is a unit } \}.
\]
is stable under multiplication and it is determined by the pair
$(\tilde{S}, \mathcal{H})$. The semigroup:
\[ \Gamma' (S):= \{ (m_1(h) , \dots, m_{c'}(h) ) \mid  h \in
\mathcal{D} \},
\]
is an analytic invariant of the germ $(S,0)$. He proved that the
semigroup $\Gamma' (S)$ is isomorphic to the projection
$\Gamma_\p'$ of the semigroup $\Gamma$ by the map
\[ \Gamma_\p \ni (a_1, \dots ,a_d)
\mapsto (a_1, \dots, a_{c'}) \in \Gamma_\p'. \] The information
encoded by the semigroup $\Gamma' (S)$, completed in certain cases
by the equisingularity type of certain transversal plane sections
to the singular locus of $S$ (a technique which is also used in
\cite{Gau}), is enough to recover a normalized sequence of
characteristic monomials of some quasi-ordinary projection of
$(S,0)$.

- The first author and Gonzalez-Sprinberg gave also a completely
algebraic proof of the analytic invariance of the sequence of
normalized characteristic monomials, which passes
 through the
invariance of the semigroup $\Gamma_\p$ (see
\cite{Pedro-Gerardo}). The approach is to use valuations defined
in an intrinsic way from the quasi-ordinary germ $(S,0)$. Those
considered in that paper are the divisorial valuations associated
to the irreducible components of the exceptional divisor obtained
after blowing up the origin $0 \in \bar{S}$, where $\bar{S}$
denotes the normalization of $S$.
 Any of these valuations induces a filtration
of $\mathcal{O}_S$ with associated graded ring isomorphic to
$\C[\Gamma_\p ]$, with a suitable graduation (notice here that the
toric variety  $Z$ mentioned before is $Z = \mbox{\rm Spec}
\C[\Gamma_\p]$). These data determines the semigroup $\Gamma_\p$
up to isomorphism by a theorem of Gubeladze (see \cite{Gu}).

Both papers \cite{patrick} and \cite{Pedro-Gerardo} used that the
normalization $\bar{S}$ of a quasi-ordinary singularity is a germ
of normal toric singularity explicitly determined by the
characteristic monomials (see \cite{Pedro-resolucion-torica} and
also \cite{patrick}). The comprehension of the phenomena however
was not completely satisfactory and unified, for instance it is
not clear why in certain cases part of the topological information
is missed in the combinatorial analysis of the preimage of the
singular locus of $S$ in $S'$. On the other hand in
\cite{Pedro-Gerardo} the use of the theorem of Gubeladze, though
useful to prove the result, does not provide a deep insight on the
reasons of the invariance.

\subsection{The zeta function of geometric monodromy}

Another motivation of this paper comes from the study of the zeta
function of geometric monodromy of the germ $f: (\C^{d+1}, 0)
\rightarrow (\C,0)$, for $f=0$ a quasi-ordinary hypersurface,
which was done by McEwan and Némethi in \cite{Nemethi} in the
irreducible case and  by these authors and González Pérez in
\cite{NemethiII} in the general case. It was proven that  the zeta
function $\z(t)$ coincides with zeta function of a plane section
$f_{x_2 = \dots = x_d =0} : ( \C^2, 0) \rightarrow (\C,0)$, if the
coordinates are suitably chosen. This result was used by Artal,
Cassou-Noguès, Luengo, and Melle in the proof of the  Monodromy
Conjecture for the local Igusa and motivic zeta functions of a
quasi-ordinary polynomial in arbitrary dimension (see
\cite{ACLM}).

The zeta function of geometric monodromy is not a complete
invariant of the embedded topological type of the germ $f=0$ if $d
> 1$. Since quasi-ordinary hypersurface singularities generalize
analytically irreducible plane curve singularities, for which the
zeta function and the Poincaré series of its semigroup coincide,
Némethi and McEwan state the following question in
\cite{Nemethi-Pregunta}:

\begin{quote}{\em "What is the right analog for $d\geq 2$
of the Gusein-Zade, Delgado and Campillo result \cite{GDC}?" }
\end{quote}

This question is also an important motivation for this paper.

\subsection{Structure of the paper and results}

We associate a  Poincaré series to a finite sequence of monomial
valuations $\underline{\nu} = (\nu_1, \dots, \nu_r)$, of a germ of
complex analytic (or algebroid) singularity $(S,0)$, when at least
one of these valuations is centered at the origin and

 - the germ  $(S,0)$ is a toric singularity of dimension $d$, defined by
$Z^{\Lambda} = \mbox{\rm Spec} \C[\Lambda]$, for $\Lambda \subset
\Z^d$ a finite type semigroup, at its special point, or,

- the germ $(S, 0)$ is an irreducible germ of quasi-ordinary
hypersurface singularity of dimension $d$.

Both types of singularities are non necessarily normal and rarely
isolated.

In our context we can not use Campillo, Delgado and Gusein-Zade
multi-index filtrations to define the Poincaré series by
(\ref{def-C-D-G}), as in \cite{Ann-toricas}, since the vector
space dimensions (\ref{cdim})  may become infinite, if one of the
monomial valuations is not centered at the origin. We deal with
this difficulty by showing that monomial valuations are {\em
compatible}, a notion which allow us to define a multi-graded
ring,
\[ \mbox{gr}_{\underline{\nu}} \mathcal{O}_S =
\bigoplus_{\underline{a} \in \Z^r_{\geq 0}} J(  \underline{a} ),
\]
associated to the analytic algebra $\mathcal{O}_S$ of the germ by
the sequence of monomial valuations $\underline{\nu} = ( \nu_1,
\dots, \nu_r)$, when at least one of them is centered at the
origin. This condition  also guarantees that the homogeneous
component $J( \underline{a} )$
 of this multi-graded
ring are vector spaces of finite dimension over $\C$. This is done
carefully in section \ref{grad-val}.  Some of the methods
introduced here might be of application to other classes of non
isolated singularities. It should be noticed that the graded ring
associated to an analytic algebra by the filtration induced by a
divisorial valuation is not necessarily Noetherian (see examples
in \cite{CS} and \cite{CGP}).

Then we define the Poincaré series of $(S,0)$ relative to the
sequence $\underline{\nu}$ as:
\begin{equation} \label{pp}
P^{\underline{\nu}} _{(S,0)} := \sum_{ \underline{a} \in
\Z^r_{\geq 0} } \dim_\C (J(  \underline{a} ) ) \,  \underline{t}^{
\underline{a} } \in \Z[[ \underline{t}]],
\end{equation}
for $\underline{t} = (t_1, \dots, t_r)$. This can be expressed as
the Poincaré series of an {\em extended semigroup} associated with
the multi-graded ring $  \mbox{gr}_{\underline{\nu}} \mathcal{O}_S
$. The Poincaré series is a rational function in $\underline{t}$
with integer coefficients.
 We exhibit then the
Poincaré series as an integral with respect to the Euler
characteristic, of a function defined by $\underline{\nu}$, on
 the projectivization of the analytic
algebra of the singularity (see Theorem \ref{euler-t} and
\ref{euler-q}).

 In the toric case,  our results generalize those obtained by  Lemahieu
 \cite{Ann-toricas}. In particular,
 if we assume that
all the monomial valuations are centered at the origin our notion
of Poincaré series,  coincides a fortiori with the one introduced
by Lemahieu \cite{Ann-toricas} using a different definition based
on (\ref{def-C-D-G}).

The class of singularities we are dealing with, toric and
hypersurface quasi-ordinary, possesses essential divisors which
are invariant by the torus action of the normalization. These
divisors define then essential divisorial valuations, which are
monomial valuations. This characterization in different cases is
done  in \cite{Bouvier}, \cite{IK}, \cite{Ishii-crelle},
\cite{Ishii-fourier} and \cite{GP-nash}. Since in our construction
we need at least one valuation centered at the origin we combine
two different notions of essential valuation, those centered at
the origin and those centered at some component of the singular
locus of $S$.

 A detailed combinatorial
description of the essential divisors in the quasi-ordinary case
in terms of the structure of the singular locus of $(S,0)$ is
given in section \ref{divi}. Of essential use is Lipman's
characterization of the singular locus of a quasi-ordinary
hypersurface (see Theorem \ref{singular}). The information on the
normalized characteristic monomials provided by the essential
valuations is encoded in a matrix $\mathcal{M}$ which is defined
and studied in section \ref{matrix}.

If $\underline{\nu}$ is the sequence of essential divisorial
valuations,  then the associated Poincaré series, $P_{(S,0)}$, is
an analytical invariant of the singularity. In the quasi-ordinary
case we show that this Poincaré series has a unique short form as
product/quotient of terms of {\em cyclotomic form}  $1 -
\underline{t}^{\underline{a}}$.

The main result of this paper is the following (see Theorem
\ref{Pseries}).

\noindent{\bf Theorem.}  {\it The Poincaré series $P_{(S,0)}$
determines and it is determined by the sequence of normalized
characteristic monomials of $(S,0)$.}

The proof is constructive, and  passes trough the analysis of this
short form of the Poincaré series, determining from it the
dimension $d$, the number of codimension one components, the
equisingular dimension $c$ and finally the characteristic
monomials by resolving a sequence of linear systems defined by
certain minors of the matrix $\mathcal{M}$.  The proof is divided
in two cases: the case $d=2$ is  simpler and similar to the one
given in \cite{PPP02}. The case $d>2$ is treated with different
arguments depending on the number $s_2$ of essential divisors over
codimension two irreducible components of the singular locus.

In comparison with the approach in \cite{patrick} the
combinatorial structure of the preimage of the singular locus in
the orbifold $S'$, which is analyzed by an smart procedure of
successive blow ups of smooth centers, missed the singularity type
of the preimage of the components in the normalization $\bar{S}$.
In our context the essential divisors keep record of this
fundamental information.  For instance, if there is only one
irreducible component of the singular locus of $S$, which in
addition is of codimension $2$, and there is  only one essential
divisorial valuation over this component then  $(S,0)$ must be
isomorphic to the germ defined by $y^2 - x_1 x_2 =0$ in $\C^d$.

We find out that the Poincaré series characterizes equisingularity
in an equisingular family of quasi-ordinary hypersurfaces (see
Corollary \ref{equi}).

In section \ref{zeta} we compare the zeta function of geometric
monodromy with the Poincaré series in the quasi-ordinary
hypersurface case and we find out that the zeta function appears
as a product of cyclotomic terms in the numerator or denominator
of the rational function in one variable obtained by specializing
the indeterminates $t_1, \dots, t_p$ in one variable $t$.
Corollary \ref{fzeta} gives a partial answer to McEwan and Némethi
question H in \cite{Nemethi-Pregunta}, and involves a new related
question.

We finish the article with an example.

\section{A reminder of toric geometry } \label{toric}

We give some definitions and notations (see \cite{Oda}, \cite{Ew}
or  \cite{Fulton} for proofs).  If $N$ is a lattice we denote by
$M$ the dual lattice, by $N_\R$ the real vector space spanned by
$N$ and by $\langle , \rangle$ the canonical pairing between the
dual lattices $N$ and $M$ (resp. vector spaces $N_\R$ and $M_\R$).
A {\it rational convex polyhedral cone} $\tau $ in $N_\R$, a {\it
  cone} in what follows, is the set $\t := \mbox{\rm pos} (a_1,
\dots, a_s)$ of non negative linear combinations of vectors $a_1,
\dots, a_s \in N$. The cone $\tau$ is {\it strictly convex} if
$\tau $ contains no linear subspace of dimension $>0$. We denote
by $\stackrel{\circ}{\tau}$ the { \it relative interior} of a cone
$\sigma$. The {\it dual} cone $\tau^\vee$ (resp. {\it orthogonal}
cone $\tau^\bot$) of $\tau$ is the set $ \{ w  \in M_\R / \langle
w, u \rangle \geq 0\}$ (resp. $ \langle w, u \rangle = 0$)  $ \;
\forall u \in \tau \}$). A {\it fan\index{fan}} $\Sigma$ is a
family of strictly convex
  cones  in $N_\R$
such that any face of such a cone is in the family and the
intersection of any two of them is a face of each. If  $\tau$ is a
cone in the fan $\Sigma$, the semigroup $\tau^\vee  \cap M$ is of
finite type, it spans the lattice $M$ and the variety $Z_{\t, N} =
\makebox{Spec} \, k [ \tau^\vee \cap M ]$, which we denote by
$Z_\tau$ when the lattice is clear from the context, is normal.
The affine varieties $Z_\tau$ corresponding to cones in a fan
$\Sigma$ glue up to define the {\it normal toric
variety\index{toric variety}} $ Z_\Sigma$. The torus $ T_N := Z_{
\{ 0 \} } \cong (k^*)^{ \mbox{rk}\, N}  $ is embedded in $
Z_\Sigma$ as an open dense subset and there is an action of $T_N$
on $Z_\Sigma$ which extends the action of the torus on itself by
multiplication.  We have a bijection between the relative
interiors of the cones  of the fan and the orbits of the torus
action, $\stackrel{\circ}{\tau} \mapsto \O_{Z_\Sigma} (\tau)  $,
which inverses inclusions of the closures. We denote the orbit
$\O_{Z_\Sigma} (\tau)$ also by $\O (\tau)$ when the toric variety
$Z_\Sigma$ is clear from the context.

Let $ \Lambda $ be a sub-semigroup of finite type of a rank $d$
lattice $ M:=\Lambda + (- \Lambda)$ which it generates as a group.
Any affine toric
  variety is of the form
$Z^{\Lambda} = \makebox{Spec} \, \C [ \Lambda ]$ where
$\C[\Lambda] = \{ \sum_{finite} c_\l X^\l \mid c_\l \in \C \}$ is
the $\C$-algebra of the semigroup. Fixing a basis of the lattice
$M$ induces an embedding of $\Lambda$ in $\Z^{d}$, corresponding
to an inclusion of $\C[\Lambda]$ in the ring of Laurent
polynomials $\C[y_1^{\pm 1}, \dots, y_d^{\pm 1} ]$, $X^{e} \mapsto
y_1^{e_1} \cdots y_d^{e_d}$, for $ (e_1, \dots, e_d)$ the
coordinates with respect to the fixed basis of the vector $e \in
\Lambda$. The torus $Z^{M}$ is an open dense subset of
$Z^\Lambda$, which acts on $Z^\Lambda$ and the action extends the
action of the torus on itself by multiplication. The semigroup
$\Lambda$ spans the cone $\s=\R_{\geq 0} \Lambda$ thus we have  an
inclusion of semigroups $\Lambda \rightarrow \bar{\Lambda} :=
\s^\vee \cap M$, defining an associated toric modification
$Z^{\bar{\Lambda}} = Z_{\s, N} \rightarrow Z^{{\Lambda}}$, which
is the {\it normalization map} (the lattice $N$ is the dual latice
of $M$) . The cone $\s^\vee$ has a vertex if and only if there
exists a zero dimensional orbit, and in this case this orbit is
reduced to the point of $ 0 \in Z^{\Lambda}$ defined by the
maximal ideal $ \mathfrak{m}_\Lambda := (\Lambda - \{ 0 \} ) \C [
\Lambda ]$. The ring $\C [[ \Lambda ]] $ is the completion of the
local ring of germs of holomorphic functions at $(Z^\Lambda, 0)$
with respect to its maximal ideal. Any non zero vector $v \in \s
\cap N$ defines a valuation $\nu$ of the field of fractions of
$\C[[\Lambda]]$, called {\em monomial
  valuation} or {\em toric valuation}, which is
defined for an element $0 \ne \f = \sum c_u X^u \in \C [[ \Lambda
]]$ by $\nu (\f) = \min_{c_u \ne 0} \langle n, u \rangle$. If the
vector $v$ is primitive the valuation $\nu$ is a  divisorial
valuation, associated to the divisor $D_v$ corresponding to the
ray $v \R_{\geq 0}$ in any suitable subdivision $\Sigma$ of $\s$,
containing this ray. Then the composite of $\p_\Sigma$ with the
normalization map of $Z^{\Lambda}$ defines a toric modification in
which the divisor $D_v$ appears and we have that the divisorial
valuation associated to $D_v$ is equal to $\nu$ (see
\cite{Bouvier}, \cite{B-GS} or \cite{Pedro-Gerardo}).

\section{Quasi-Ordinary Singularities}\label{quasi-ord-sect}

A germ of algebroid hypersurface  $(S,0) \subset (\C^{d+1}, 0) $
is {\it quasi-ordinary} if there exists a finite morphism $(S, 0)
\rightarrow (\C^d,0)$ (called a {\it quasi-ordinary
projection\index{quasi-ordinary projection}}) such that the
discriminant locus is contained (germ-wise) in a normal crossing
divisor. In suitable coordinates depending on this projection, the
hypersurface $(S, 0)$ has an equation $f=0$, where $f \in \C [[
\underline{x} ]][y]$ is a {\it quasi-ordinary polynomial}: a
Weierstrass polynomial with discriminant, $\D_y f$, as a
polynomial in $y$, of the form $\D_y f = \underline{x}^\d
\epsilon$, where $\epsilon$ is a unit in the ring $ \C [[
\underline{x} ]]$ of formal power series  in the variables
$\underline{x}= (x_1, \dots, x_d)$ and  $\d \in \Z^d_{\geq 0}$.

\subsection{Characteristic data associated to a quasi-ordinary
projection}\label{qo}

We will suppose from now on that the germ $(S, 0)$ is irreducible,
i.e., the polynomial $f$ is irreducible. In this subsection we
recall some numerical data which are associated to the fixed
quasi-ordinary projection, or equivalently to the quasi-ordinary
polynomial $f$. We use these data to describe the singular locus
and some invariants of the germ $(S,0)$.

The Jung-Abhyankar theorem guarantees that all the roots $\{ \z
^{(l)} \}_{l=1}^n $ of $f$, which are called {\em quasi-ordinary
branches}, are fractional power series in the ring $\C \{
\underline{x}^{1/n} \}$ for $\underline{x}^{1/n} = (x_1^{1/n},
\dots, x_d^{1/n}) $ and $n := \deg f \in \Z_{\geq 0}$ (see
\cite{Abhyankar}). Since the discriminant of $f$ is equal to $\D_y
f = \prod_{i\ne j} (\z^{(i)} - \z^{(j)})$, each factor $\z^{(i)} -
\z^{(j)}$ is of the form a monomial times a unit in $\C [[
\underline{x}^{1/n} ]]$. These monomials (resp. their exponents)
are called {\it characteristic} or {\it distinguished}, of the
quasi-ordinary polynomial $f$, i.e., these exponents depend on the
quasi-ordinary projection considered.

The characteristic exponents can be relabelled in the form
\begin{equation} \label{ch-order}
\l_1 < \l_2 < \cdots < \l_g
\end{equation}
where $<$ means $\leq$  coordinate-wise and $\ne$  (see
\cite{Lipman88}).

Without loss of generality we relabel the variables $x_1, \dots,
x_d$ in order to have that if,  $\l_j = ( \l_j^1, \dots, \l_j^d)
\in \Q^d$ for $j=1, \dots, g$,
\begin{equation} \label{lex}
(\l_1^1, \dots, \l_g^1) \geq_{\mbox{\rm lex}} \cdots
\geq_{\mbox{\rm lex}} (\l_1^d, \dots, \l_g^d),
\end{equation}
where $\mbox{\rm lex}$ is lexicographic order.

\begin{definition} \label{lattices}
The characteristic exponents determine the following nested
sequence of lattices:
 \[  M_0 := \Z^d \subset M_1 \subset \cdots \subset M_g =: M \]
where $M_j := M_{j-1} + \Z \l_j$ for $j=1, \dots, g$ with the
convention $\l_{g +1} = + \infty$. The index $n_j$ of the lattice
$M_{j-1}$ in the lattice $M_j$ is $>1$. We call $n_1, \dots, n_g$
the {\em characteristic integers} associated to $f$ (see
\cite{Lipman88}, \cite{Pedro-Tesis}, \cite{Pedro-semigrupo} and
\cite{patrick}).
\end{definition}

Since the polynomial $f$ is irreducible, we can identify the
analytic algebra $\mathcal{O}_S:= \C[[ \underline{x}]] [y]/(f)$ of
the germ $(S, 0)$ with the ring $ \C [[ \underline{x} ]] [ \z] $,
for $ \z$ any fixed quasi-ordinary branch $ \z$ parameterizing
$(S, 0)$.

We denote by $N$ (resp. by $N_j$) the dual lattice of $M$ (resp.
of $M_j$, for $j=1, \dots,g$) and by   $\s \subset N_\R$ the cone
spanned by the dual basis of  the canonical basis of $M_0$. With
these notations we have that the homomorphism $\C[[\underline{x}]]
\rightarrow \C[[\s^\vee \cap M_0 ]] $ which maps $x_i \mapsto
X^{e_i} $, for $e_i$ running through the primitive integral
vectors of the lattice $M_0$ in the edges of $\s^\vee$, is an
isomorphism of $\C$-algebras since $\{ e_i \}_{i=1, \dots, d}$ is
a basis of $M_0$. Remark that $X^{e}$ denotes  a monomial of the
semigroup algebra, without making any reference to a choice of
basis of the lattice. Since the exponents of the quasi-ordinary
branch $\z$ belong to the semigroup $\s^\vee \cap M$ we obtain a
ring extension $ \mathcal{O}_S = \C [[ \s^\vee \cap M_0 ]][\z]
\rightarrow \C [[ \s^\vee \cap M ]]$, which is the inclusion of
$\mathcal{O}_S $ in its integral closure in its field of
fractions.
\begin{proposition}\label{normalization} (see Proposition
14, \cite{Pedro-resolucion-torica}) The map $n: (Z_{\s, N}, o_\s)
\rightarrow (S,0)$, corresponding to the ring extension
\begin{equation} \label{extension}
\C [[ \s^\vee \cap M_0 ]] [\z] \rightarrow \C [[ \s^\vee \cap M
]],
\end{equation}
is the normalization map.
\end{proposition}
More generally the normalization of a germ of quasi-ordinary
singularity, non necessarily hypersurface, is a toric singularity
(see \cite{patrick}).

The elements of $M$ defined by:
\begin{equation}\label{rel-semi}
{\g}_1 =  \l_1 \quad \mbox{ and } \quad {\g}_{j+1} - n_j {\g}_{j}
= \l_{j+1} -  \l_{j}, \quad \mbox{ for } \quad j= 1, \dots, g-1,
\end{equation}
together with the semigroup $\Gamma_\p := \Z^d_{\geq 0} +  \g_1
\Z_{\geq 0} + \cdots + \g_g \Z_{\geq 0} \subset \s^\vee \cap M_\Q$
are associated with the quasi-ordinary polynomial $f$ (see the
introduction for more on the notion of semigroup in this case).

\begin{remark} \label{function}
For any $\g \in \Gamma_\p$ there exists a function $\f_\g \in
\mathcal{O}_S $ such that $\f_\g = X^{\g} \cdot \epsilon_\g$ in
$\C[[\s^{\vee} \cap M]]$, where $\epsilon_\g$ is a unit in
$\C[[\s^{\vee} \cap M]]$ (see \cite{Pedro-semigrupo},
\cite{Pedro-resolucion-torica} or \cite{patrick}).
\end{remark}

We have that for $1 \leq j < g$
\begin{equation} \label{preforma}
n_j\gamma_{j} < \g_{j+1},
\end{equation}
where $<$ means $\ne$ and $\leq$ coordinate-wise. The equality
\begin{equation} \label{forma}
\g = \a  + l_1  \g_1 + \cdots + l_g  \g_g, \quad \mbox{ with }
\alpha  \in \s^\vee \cap M_0, \, 0\leq l_i  < n_i, \, i = 1,
\dots, g;
\end{equation}
determines an element $\g$ of the semigroup $\Gamma_\p$ and any
element of this semigroup has a unique expansion of the form
(\ref{forma}), (see \cite{Pedro-semigrupo}  Lemma 3.3).

\begin{definition}
The {\em generating series of the semigroup $\Gamma_\p$} is the
series
\[
P_{\Gamma_\p} := \sum_{\g \in \Gamma_\p} X^{\gamma} \in
\C[[\Gamma_\p]].
\]
\end{definition}

\begin{proposition} \label{naive}
The generating series of the semigroup $\Gamma_\p$ is a rational
function:
\[
P_{\Gamma_\p}= \prod_{j=1}^d \frac{1}{1- X^{e_j}} \prod_{i=1}^g
\frac{1- X^{n_i \g_i} }{1 - X^{\gamma_i}} \in \C(\Gamma_\p),
\]
where $\C(\Gamma_\p)$ denotes the field of fractions of
$\C[\Gamma_\p]$.
\end{proposition}
\textbf{Proof:} We use the unique expansion form (\ref{forma}) of
the elements of the semigroup:
\[
\begin{array}{c}
P_{\Gamma_\p}  =   \sum_\g X^{\a (\g)} X^{\sum l_i (\g) \g_i}
 \stackrel{\mbox{\small(\ref{forma})}}{=} \sum_{\a \in \s^\vee \cap
 M_0} X^\a \prod_{i=1}^g ( 1 + X^{\g_i} + \dots + X^{(n_i-1)\g_i}).
\end{array}
\]
Notice by construction the semigroup $\s^\vee \cap M_0$ is
isomorphic to $\Z^d_{\geq 0}$. Then, rewriting this expression
using elementary properties of the geometric series gives the
result. $\Box$

In the plane branch case we have a similar computation (see
\cite{C-D-G-curvas}).

\begin{definition}
The characteristic exponents $\l_1, \dots, \l_g$ of a
quasi-ordinary polynomial $f$ are {\it normalized} if (\ref{lex})
is verified, and in addition if $\l_1= (\l_1^1, 0, \dots, 0)$ then $\l_1^1
>1$.
\end{definition}
If this is the case we say that the quasi-ordinary polynomial $f$
or the associated quasi-ordinary projections are also {\em
normalized}. Lipman showed how to associate a normalized
quasi-ordinary polynomial $f'$, to a given quasi-ordinary
polynomial $f$, in such a way that $f=0$ and $f'=0$ define the
same germ and the characteristic exponents associated to $f'$, are
related with those associated to $f$ by certain inversion formulas
(see \cite{Gau} Appendix or \cite{Lipman65}).

\begin{remark} We assume that the characteristic
exponents and quasi-ordinary polynomials are normalized unless the
contrary is stated.
\end{remark}

Suppose that the reduced discriminant locus of the quasi-ordinary
projection $\p : (S,0) \rightarrow (\C^d,0)$ has equation $x_1
\cdots x_{c} =0$. The following result is a reformulation of part
of Theorem 7.3, \cite{Lipman88}, using our notations.

\begin{theorem} \label{singular}
With the previous notations, we have that the singular locus
$\mbox{\rm Sing } S$ of the quasi-ordinary hypersurface $(S,0)$
has only irreducible components of codimension one or two of the
form $Z_i := S \cap \{x_i =0 \} $ or $Z_{i,j} := S \cap \{ x_i =
x_j =0 \}$,  for $1\leq i, j \leq c$. If $1 \leq i \leq c$, the
coordinate section $Z_i$ is a component unless $\l_1^i = \dots =
\l_{g-1}^i = 0 $ and $\l_g^i = 1/ n_g$. If $1\leq i, j \leq c$ ,
the coordinate section $Z_{i,j}$ is a component if and only if
neither $Z_i$ nor $Z_j$ is a component of the singular locus of
$S$.
\end{theorem}

\section{The graded ring associated with finitely many
monomial valuations} \label{grad-val}

We define a graded ring associated to the analytic algebra of a
toric or a quasi-ordinary singularity by a sequence of toric
divisorial valuations. We continue in this section to extend the
approach and the results obtained in \cite{Pedro-Gerardo}
(corresponding to the case of a single valuation).

A discrete valuation $\nu$ of the field of fractions of a ring
$R$, which is centered on $R$, defines a filtration of $R$ with
ideals $I(b) =  \{ \f \in R \mid \nu ( \f) \geq b \}$, for $b \in
\Z_{\geq 0}$. The associated graded ring is $\mbox{gr}_\nu R :=
\bigoplus_{b \geq 0} I(b) / I(b+1)$.

Suppose now that the ring $R$ is graded over the semigroup
$\Z^{r}_{\geq 0}$, i.e., it is of the form:
\[
R = \bigoplus_{\underline{a} \in \Z^r_{\geq 0}} J(\underline{a})
\]
where $J(\underline{a})$ denotes the homogeneous component
associated to $\underline{a} \in \Z_{\geq 0}^r$. In this case we
say that $R$ is {\em multi-graded} or {\em simply} that $R$ is
graded. In general an ideal $I(b) \subset R$ associated to the
valuation is not homogeneous with respect to the grading of $R$.

\begin{definition}
We say that the valuation $\nu$ is {\em compatible} with the
graduation if for all   $b \geq 0$ the ideal $I(b)$ is an
homogeneous ideal, i.e.,
\[
I(b) = \bigoplus_{\underline{a} \in \Z^r_{\geq 0}}
I(\underline{a},b), \mbox{ where } I(\underline{a},b) := \{ \f \in
J(\underline{a}) \mid \nu(\f) \geq b \}.
\]
\end{definition}
It follows that  the associated graded ring $\mbox{gr}_\nu R $
inherits a multi-grading since
\[
I(b)/I(b+1) = \bigoplus_{a \in \Z^r_{\geq 0 }}  I(\underline{a},b)/
 I(\underline{a}, b+1).
\]
If we denote by $J(\underline{a},b) := I(\underline{a},b)/
 I(\underline{a}, b+1)$ then we have that:
\[
 \mbox{gr}_\nu R := \bigoplus_{\underline{a} \in \Z^r_{\geq 0} ,b \geq 0}
 J (\underline{a},b).
 \]

\begin{remark}
This notion of compatibility of a valuation with a graded ring can
be defined more generally for non discrete valuations and for
gradings on different semi-groups, though we do not need this in
this paper.
\end{remark}

\subsection{Toric case}

Let $(Z^\Lambda,0)$ be the germ of an affine toric variety at its
zero orbit, which is assumed to exist. The normalization
$Z^{\bar{\Lambda}} $ is of the form $Z^{\bar{\Lambda}} = Z_{\s, N}
$.

Let  $ \underline{v} =(v_1, \dots, v_r)$ be a finite sequence of
non zero vectors in $\s \cap N$, we consider then  the {\em
graduation} of the $\C$-algebra $\C[ \Lambda ]$,
\begin{equation} \label{graded}
\C [ \Lambda ] ^{(\underline{v}) } := \bigoplus_{\underline{a} \in
\Z^r_{\geq 0}} H (\underline{a})
\end{equation}
 where the  homogeneous terms $H (\underline{a})$ are
 defined by:
\[
H ( \underline{a} ) := \bigoplus^{u \in \Lambda}_{ \langle v_j, u
\rangle = a_j, j=1, \dots, r} \C X^u, \quad \mbox{ for }
\underline{a} = (a_1, \dots, a_r) \in \Z^r_{ \geq 0}.
\]

We show that the graduation $ \C [ \Lambda ] ^{(\underline{v} ) }
$ introduced above is in fact a multi-grading associated to the
sequence of {\em monomial valuations} $\nu_1, \dots, \nu_r$
corresponding to a the sequence of vectors $v_1, \dots, v_r \in \s
\cap N - \{ 0\}$.

Let $v \in \s \cap N$ be a non zero vector defining a monomial
valuation $\nu$ of the ring $R$, equal to $\C[[\Lambda]]$ (or to
$\C[\Lambda]$). The associated graded ring is  $\mbox{gr}_\nu R =
\bigoplus_{a \geq 0}  J(a)$ where $J(a) = I(a) / I(a+1)$ and
$I(a)$ denotes the ideals of the filtration associated to $\nu$.
If $0 \ne h \in \C [[\Lambda]]$ then $h_{|
  v}$
denotes the {\em symbolic restriction} of $h$ to the face defined
by the vector $v$ on its Newton polyhedron (see
\cite{Pedro-Gerardo}).

\begin{proposition}  \label{key} With the above notations the mapping:
\[
J(a) - \{ 0 \} \longrightarrow H(a) - \{ 0 \}, \quad h \mbox{ mod
} I(a+1) \mapsto h_{| v},
\]
extends to an isomorphism of graded rings:
\begin{equation} \label{iso-gr}
\psi: \mbox{gr}_\nu R \longrightarrow \C[\Lambda]^{(v)}.
\end{equation}
If $0 \ne w \in \s \cap N$ defines a monomial valuation $\omega$
of $R$ then define \[ \bar{\omega} : \mbox{gr}_\nu R - \{ 0 \}
\longrightarrow \Z_{\geq 0} \] on homogeneous elements by,
\[
J(a) - \{ 0 \} \longrightarrow \Z_{\geq 0}, \quad \bar{h} \mapsto
\max \{ \omega (h) \mid h + I(a+1) = \bar{h} \}
\]
and in general if $0 \ne \f= \sum_a  \bar{h}_a \in \mbox{gr}_\nu R
$ is a sum of homogeneous terms $0 \ne \bar{h}_a \in J(a)$ then
define
\[
\bar{\omega}(\f) := \min_a \{ \bar{\omega} (\bar{h}_a)\}.
\]
Then we have that $\bar{\omega} = \omega \circ \psi$, hence
$\bar{\omega}$ is a monomial valuation of $\mbox{gr}_\nu R$
 compatible with the graded
structure.
\end{proposition}
\textbf{Proof:} For the first assertion see the proof of
\cite{Pedro-Gerardo} Proposition 3.2. For the second assertion we
use that $\psi$ is an isomorphism of graded rings, hence induces
an isomorphism of homogeneous components: $J(a) \rightarrow H(a)$.

If $0 \ne \bar{h} \in  J(a)$ let us take a representative $h_0 \in
I(a)- I(a+1)$, then  we have $\psi (\bar{h}) = (h_0 )_{| v}  \in
H(a)$. Notice also that $ (h_0 )_{| v}  + I(a+1) = \bar{h}$ and
that $\max \{ \w ( (h_0 )_{| v} +  h') \mid h ' \in I(a+1)  \} $
is reached, for instance, when $h' =0$. This implies that on each
homogeneous component $J(a)$ we have that $\bar{\omega} = \omega
\circ \psi$. Notice now that both maps $\omega \circ \psi$ and
$\bar{\omega}$ coincide on $\mbox{gr}_\nu R $. This implies that
$\bar{\omega}$ is a monomial valuation of $\mbox{gr}_\nu R $, with
respect to the toric structure defined by $\psi$. $\Box$

We call $\bar{w}$ an {\em induced valuation} on the graded ring
$\mbox{gr}_\nu R $.

\begin{proposition} \label{key3}
Let $\underline{\nu}= (\nu_1, \dots, \nu_r)$ be a sequence of
monomial valuations of the ring $R$, equal to $\C[\Lambda] $ or
$\C[[\Lambda]]$, then these valuations induced a multi-graded ring
$\mbox{gr}_{\underline \nu} R $ given with a sequence of
valuations $\underline{\bar{\nu}} = ( \bar{\nu}_1, \dots,
\bar{\nu}_r) $ on it, and there   exists an isomorphism of graded
rings
\[
\psi : \mbox{gr}_{\underline \nu} R  \rightarrow \C
[\Lambda]^{(\underline{v})},
\]
such that $\nu_i \circ \psi = \bar{\nu}_i$ for $i=1, \dots, r$.
\end{proposition}
\textbf{Proof:} We generalize the proof of Proposition \ref{key}
by induction on $r$. At step $k$ we have by induction hypothesis a
graded ring $\mbox{gr}_{\underline (\nu_1, \dots, \nu_k)} R$ with
valuations $\tilde{\nu}_i$, for $i=1, \dots, r$ and a graded
isomorphism
\[
\tilde{\psi}: \mbox{gr}_{ (\nu_1, \dots, \nu_k)} R \rightarrow \C
[\Lambda]^{(v_1, \dots, v_k)},
\]
such that  $\nu_i \circ \tilde{\psi} = \tilde{\nu}_i$, for $i=1,
\dots, r$. This implies that $\tilde{\nu}_i$ can be seen as a
monomial valuation of $\mbox{gr}_{ (\nu_1, \dots, \nu_k)} R$ which
is compatible with the graded structure. If $k <r$ it follows from
Proposition \ref{key} that $\tilde{\psi}$ induces a graded
isomorphism:
\[
\bar{\psi}: \mbox{gr}_{ (\nu_1, \dots, \nu_{k+1})} R \rightarrow
\C [\Lambda]^{(v_1, \dots, v_{k+1})},
\]
which is given in homogeneous components, with the above notations
by,
\[
{J} (a_1, \dots, a_{k+1}) - \{ 0 \}  \longrightarrow H (a_1,
\dots, a_{k+1})  - \{ 0 \},  \quad h \ \mbox{mod} \ I_{\nu_{k+1} }
(a_{k+1} +1)  \mapsto (\tilde{\psi} (h))_{| v_{k+1}}.
\]
Notice that this is a natural way of extending  the map $\psi$
defined in Proposition \ref{key}, and makes of $\bar{\psi}$ an
isomorphism of the graded rings quoted above. Then the last
assertion on the induced valuations is implied by Proposition
\ref{key}. $\Box$

\begin{remark}
Taking a different order on the monomial valuations in Proposition
\ref{key3} provides isomorphic pairs of graded rings with induced
valuations. The graded ring $\mbox{gr}_{\underline{\nu}} (R)$ is
graded on the semigroup \[ \langle \underline{v}, \Lambda
\rangle = \{ \langle \underline{v}, \g \rangle :=  (\langle
\underline{v}_1, \g \rangle, \dots, \langle \underline{v}_r,
\g \rangle )  \mid \mbox{ for } \g \in \Lambda \} \subset
\Z^r_{\geq 0}.
\]
\end{remark}
Indeed, we have that $\underline{a} \in \langle \underline{v},
\Lambda \rangle$ if and only if $J(\underline{a} ) \ne (0)$.

\subsection{Quasi-ordinary hypersurface case}

Let $(S,0)$ be an irreducible germ of quasi-ordinary hypersurface
given with a normalized quasi-ordinary projection $\p$, inducing a
toric structure on the normalization $(\bar{S},0) = ( Z_{\s, N},
o_\s)$. Any vector $v \in \s \cap N$ defines a {\em monomial
valuation} on $\C[[\s^\vee \cap M]]$, and hence on  $\mathcal{O}_S
\subset \C[[\s^\vee \cap M]] $ by the inclusion (\ref{extension}).
The associated graded ring is  $\mbox{gr}_\nu \mathcal{O}_S =
\bigoplus_{a \geq 0}  J(a)$ where $J(a) = I(a) / I(a+1)$ and
$I(a)$ denotes the ideals of the filtration associated to $\nu$.

\begin{proposition}\label{key-bis}
Let $\underline{\nu}= (\nu_1, \dots, \nu_r)$ be a sequence of
monomial valuations of the ring $R= \mathcal{O}_S $, with at least
one, say $\nu_1$ centered at the maximal ideal of $R$. This
sequence of valuations induce a multi-graded ring
$\mbox{gr}_{\underline \nu} R $ given with a sequence of
valuations $\underline{\bar{\nu}} = ( \bar{\nu}_1, \dots,
\bar{\nu}_r) $ on it, and  there exists an isomorphism of graded
rings
\[
\psi : \mbox{gr}_{\underline \nu} R  \rightarrow \C
[\Gamma_\p]^{(\underline{v})},
\]
such that $\nu_i \circ \psi = \bar{\nu}_i$ for $i=1, \dots, r$.
\end{proposition}
\textbf{Proof:} The sequence $\underline{\nu}$ corresponds to a
sequence of vectors $\underline{v} =(v_1, \dots, v_r)$ in $\s \cap
N$, with $v_1 \in \stackrel{\circ}{\s} \cap N$. The assertion
corresponding to the case of one monomial valuation, necessarily
centered at the maximal ideal of $R$ is similar to the statement
of Proposition \ref{key}. For the proof see Proposition 3.2 of
\cite{Pedro-Gerardo}. Remark that the hypothesis, $v_1 \in
\stackrel{\circ}{\s} \cap N$, is essential in this case. The rest
of the assertion follows formally from the proof of Proposition
\ref{key3}. $\Box$

\begin{remark} \label{permutation}
If $\t$ is a permutation of $\{1, \dots, r\}$ such that $v_{\t(1)
} \in \stackrel{\circ}{\s} \cap N$ then the sequences of vectors
$\underline{v}$ and  $\underline{v}_{\t} := (v_{\t(1)}, \dots,
v_{\t(r)})$ define isomorphic pairs of graded rings with induced
valuations. The graded ring $\mbox{gr}_{\underline{\nu}}
(\mathcal{O}_S)$ is graded on the semigroup \[ \langle
\underline{v}, \Gamma_\p \rangle := \{ \langle \underline{v},
\g \rangle :=  (\langle \underline{v}_1, \g \rangle, \dots,
\langle \underline{v}_r, \g \rangle )  \mid \mbox{ for } \g \in
\Gamma_\p \} \subset \Z^r_{\geq 0}. \]
\end{remark}

\section{The Poincaré series as an integral with respect to the Euler
characteristic}\label{integral}
Given a germ of affine toric singularity $(Z^{\Lambda},0)$ or an
irreducible germ of quasi-ordinary singularity  $(S,0)$ provided
with a sequence of monomial valuations  $\underline{\nu} = (\nu_1,
\dots, \nu_r) $ (with at least one say $\nu_1$ centered at the
origin) we define a notion of Poincaré series from the extended
semigroup associated with the graded ring
$\mbox{gr}_{\underline{\nu}} R$, for $R = \C[[\Lambda]]$ or $R=
\mathcal{O}_S$. Then we show that this Poincaré series is an
integral with respect to the Euler characteristic over the
projectivization of $R$ of certain function defined by
$\underline{\nu}$.

\subsection{Toric case}

We consider a notion of {\em extended semigroup},
\[
\widetilde{\Lambda}: = \{ (\underline{a}, \f) \mid
\underline{a}\in \langle \underline{v}, \Lambda \rangle, \mbox{
and } \f \in J(\underline{a}) \},
\]
associated to the graded ring
\[
\mbox{gr}_{\underline{\nu}} R  = \bigoplus_{ \underline{a} \in
\langle \underline{v}, \Lambda
 \rangle} J( \underline{a}), \mbox{ for } R =\C[[\Lambda]] \mbox{
 or } \C[\Lambda],
 \]
following Campillo, Delgado and Gusein-Zade  (see
\cite{C-D-G-curvas}, \cite{GDC},
\cite{C-D-G-superficies-racionales} and also \cite{Ann-toricas}).
The operation in the semigroup $\widetilde{\Lambda}$ is additive
in the first entry and multiplicative in the second. The
projection $\r : \widetilde{\Lambda} \rightarrow \Z^r$ has fibers
$\r^{-1} (\underline{a}) = J( \underline{a})$ if $ \underline{a}
\in \langle \underline{v}, \Lambda \rangle $ and $\r^{-1}
(\underline{a}) = \{ 0 \}$ otherwise.

If $L$ is a $\C$-vector space we denote by $\P L$ its
projectivization, notice that if the dimension of $L$ is finite,
then it is equal to the Euler characteristic $\chi (\P L)$.

If the sequence of monomial valuations  $\underline{\nu} = (\nu_1,
\dots, \nu_r) $ corresponds to the sequence of vectors
$\underline{v}= (v_1, \dots, v_r)$, then our hypothesis is
equivalent to $v_1 \in \stackrel{\circ}{\s} \cap N$. This implies
that for any $\underline{a}$ we have:
\begin{equation} \label{chi}
 \dim_\C J( \underline{a}) =  \chi (\P  J( \underline{a}) ) <
\infty.
\end{equation}

\begin{definition}
The {\em Poincaré series} associated to the germ $(Z^\Lambda, 0)$
and the sequence of monomial valuations $\underline{\nu}$:
\begin{equation} \label{chi-ps}
P^{\underline{\nu}} _{ (Z^\Lambda, 0)}  := \sum_{\underline{a} \in
\Z^r_{\geq 0} }   \chi (\P  J( \underline{a}) )
\underline{t}^{\underline{a}}
\in \Z[[ \underline{t}]].
\end{equation}
\end{definition}
\begin{remark} We often denote $P^{\underline{\nu}} _{ (Z^\Lambda, 0)}$ by
$P^{\underline{v}} _{ (Z^\Lambda, 0)}$ , for $\underline{v}$ the
sequence of vectors in $\s \cap N$ defining the sequence of
monomial valuations $\underline{\nu}$.
\end{remark}
\begin{theorem} \label{caso-torico}
Let $\underline{v} = (v_1, \dots, v_r)$ be a sequence of vectors
in $\s  \cap N$, with $v_1 \in \stackrel{\circ}{\s}$, defining the
the monomial mapping $\Phi^{\underline{v}}$:
\begin{equation} \label{phi}
\Phi^{\underline{v}}: \C[[\Lambda]] \longrightarrow
\C[[\underline{t}]], \quad X^{\l} \mapsto \underline{t}^{\langle
\underline{\nu}, \l \rangle} := t_1^{\langle v_1, \l \rangle}
\cdots t_r^{\langle v_r, \l \rangle}.
\end{equation}
Then we have:
\[
\Phi^{\underline{v}} (   \sum_{\l \in \Lambda} X^{\l}     ) =
P^{\underline{v}}_{(Z^\Lambda, 0)}.
\]
\end{theorem}
\textbf{Proof:} This is a direct consequence of the existence of a
graded isomorphism between $\mbox{gr}_{\underline{v}} R$ and
$\C[\Lambda]^{\underline{v}}$ (see Proposition \ref{key3}), since
the number $ \chi (\P  J( \underline{a}) ) $ is equal to the
dimension of the $\C$-vector space $J( \underline{a})  $, which is
isomorphic to $H(\underline{a})$. It follows that this dimension
 is the number
of elements in the set $\{ \l \in \Lambda \mid \langle v_j, \l
\rangle = a_j, j=1, \dots, r \}$ (which is assumed to be finite
for every $\underline{a} \in \Z^{r}_{\geq 0}$). $\Box$

\begin{remark} The notion of Poincaré series associated to a finite
number of monomial valuations on a toric variety $Z^{\Lambda}$ was
studied by Lemahieu \cite{Ann-toricas} using only valuations
centered at $0 \in Z^\Lambda$ and the definition of multi-index
filtration of the form (\ref{def-C-D-G}). The same statement of
Theorem \ref{caso-torico} holds, with a different proof (see
Theorem 1 in \cite{Ann-toricas}). The definition of the Poincaré
series presented here (when at least one of the vectors of the
sequence $\underline{v}$ belongs to $\stackrel{\circ}{\s}$),
generalizes her definition (when all the vectors in the sequence
$\underline{v}$ belong to $\stackrel{\circ}{\s}$).
\end{remark}

\begin{remark}
It is well known that the generating series of the semigroup
$\Lambda$, $P = \sum X^{\l} \in \C[[\Lambda]] $ is a rational
function, i.e., it belongs to the ring of fractions of
$\C[\Lambda]$ (see \cite{Stanley}, apply to this particular case
Theorem 8.20 \cite{Miller},  or \cite{Barvinok} for an algorithmic
approach). It follows from Theorem \ref{caso-torico} that the
Poincaré series $P_{(Z^\Lambda, 0)} ^{\underline{\nu}}$ is also a
rational function in $\underline{t}$.
\end{remark}

Following Campillo, Delgado and Gusein-Zade, we exhibit the
Poincaré series as an integral with respect to the Euler
characteristic of the projectivization of the local ring $R =
\C[[\Gamma]]$, considered as a $\C$-vector space.

For any $k \geq 0$ we denote by
\begin{equation} \label{val-jets}
\mathcal{J}_R^k (\nu_1) : = R / I_{\nu_1} (k+1)
\end{equation} the
$k$-jet of functions in $R$, relative to the divisorial valuation
$\nu_1$ (where $I_{\nu_1} (k+1) = \{ \f \in R \mid \nu_1 (\f) \geq
k+1 \}$). This is the analogous for a smooth source of the
classical $k$-jet of functions $\mathcal{O}_{\C^d, 0} /
\mathfrak{m}^{k+1}$. The important fact is that $\nu_1$ is
centered at the origin, hence the $k$-jet space $\mathcal{J}_R^k
(\nu_1)$ is a $\C$-vector space of finite dimension $d_{\nu_1}
(k)$. We denote often    $\mathcal{J}_R^k (\nu_1)$ and $d_{\nu_1}
(k)$  by $\mathcal{J}_R^k$ and $d(k)$ respectively, for short.

The symbol $\P^* L$ denotes the disjoint union of $\P L$ and a
point $*$, for $L$ a finite dimensional vector space. We have
natural maps $\p_k : \P R \rightarrow \P^* \mathcal{J}_R^k
(\nu_1)$ and $\p_{k,k'} : \P^* \mathcal{J}_R^k (\nu_1) \rightarrow
\P^* \mathcal{J}_R^{k'} (\nu_1) $, for $k \geq k' \geq 0$.

\begin{definition} \label{cylindric}
A subset $X \subset \P R$ is {\em cylindric} if $X = \p_k ^{-1}
(Y)$, for $Y$ some constructible subset $Y  \subset \P
\mathcal{J}_R^{k} \subset  \P^* \mathcal{J}_R^{k}$. The {\em Euler
characteristic of the cylindric set $X$} is defined as $\chi(X) :=
\chi (Y)$.
\end{definition}
This notion is well defined since the map $\p_{k,k'}$ over  $ \P
\mathcal{J}_R^{k'} \subset  \P^* \mathcal{J}_R^{k'}$ is a locally
trivial fibration, with fibers complex affine spaces of dimension
$d(k) - d(k')$.

Let $\Phi: \P R \rightarrow G$ be a function with values in an
abelian group $G$. We say that $\Phi$ is cylindric if for each $g
\in G$, $g \ne 0$ the set $\Phi^{-1} (g) \subset \P R$ is
cylindric.  The integral of the cylindric function $\Phi $ over
the space $\P R$ with respect to the Euler characteristic is
defined if the sum in the right-hand side of the following formula
makes sense in $G$,
\[
\int_{\P R} \Phi d \chi := \sum_{g \in \Phi(\P R), g \ne 0} \chi (
\Phi^{-1} (g) ) g.
\]
In this case we say that the function $\Phi$ is integrable.

We apply this construction for $G= \Z[[\underline{t}]]$ and $\Phi
:= t^{\underline{\nu}}$, which is well defined on $\P R$.
\begin{theorem} \label{euler-t}
With the above notations we have that:
\[
P_{(Z^\Lambda, 0) }^{\underline{\nu}} = \int_{\P R}
t^{\underline{\nu}} d \chi.
\]
\end{theorem}
\textbf{Proof:} Notice that $\Phi (\P R) = \{ t^{\underline{a}}
\mid \underline{a} \in \langle \underline{\nu}, \Lambda \rangle
\}$. Notice also that if $\underline{a} \in \langle
\underline{\nu}, \Lambda \rangle $ then
\begin{equation} \label{asteris}
\Phi^{-1} (t^{\underline{a}}) = \bigcap_{i} I_{\nu_i} (a_i)
\setminus \bigcup_{i} I_{\nu_{i}} (a_{i} + 1).
\end{equation}
This set is cylindrical since
\[
\Phi^{-1} (t^{\underline{a}}) =\p_{a_1}^{-1} ( \P H
(\underline{a}) ), \] where  $  H (\underline{a}) $
 is seen as a subset of $\mathcal{J}_R^{a_1} (\nu_1)$.
Then it follows that
\[
\chi ( \Phi^{-1} (t^{\underline{a}})) = \chi ( \P
H(\underline{a})) = \chi ( \P J(\underline{a}) )
\]
by the isomorphism  $H(\underline{a}) \rightarrow J(\underline{a})
= \r^{-1} (\underline{a}) $ of Proposition \ref{key3}. We have
deduced that
\[
\int_{\P R} t^{\underline{\nu}} d \chi = \sum_{\underline{a} \in
\langle \underline{\nu}, \Lambda \rangle } \chi ( \P
J(\underline{a}) ) t^{\underline{a}} = P_{(Z^\Lambda, 0)
}^{\underline{\nu}} . \quad \Box
\]

\begin{remark}
A fortiori it follows that the integral $\int_{\P R}
t^{\underline{\nu}} d \chi$ is independent of the valuation
$\nu_1$ centered at the origin, chosen among the sequence
$\underline{\nu}$ to define the cylinders of $\P R$.
\end{remark}

\subsection{The Poincaré series of a quasi-ordinary hypersurface}

Let $(S,0)$ denote an irreducible  germ of quasi-ordinary
hypersurface defined by a normalized quasi-ordinary polynomial
$f$, and with associated toric normalization of the form
$(\bar{S},0) = (Z_{\s, N}, o_\s)$ (see section \ref{qo}).

Let $\underline{v} = (v_1, \dots, v_r)$ be a sequence of  vectors
in $\s \cap N$ with $v_1 \in \stackrel{\circ}{\s}$ defining
monomial valuations  $\underline{\nu} = (\nu_1, \dots, \nu_r)$ and
a multi-graded ring \[
 \mbox{\rm gr}_{\underline{v}} \mathcal{O}_S
= \bigoplus_{\underline{a} \in \langle \underline{v}, \Gamma_\p
\rangle } J(\underline{a})\] (see Proposition \ref{key-bis}).
 Since $v_1 \in \stackrel{\circ}{\s}$ we have that formula
(\ref{chi}) holds (the hypothesis $v_1 \in \stackrel{\circ}{\s}$
 is also needed in Proposition
\ref{key-bis}). We can define analogously the  {extended
semigroup} $\widetilde{\Gamma}$, and the associated {\em Poincaré
series} $P^{\underline{\nu}} _{ (S,0)} $ by formula
(\ref{chi-ps}).

By Proposition \ref{key-bis}  the graded ring
$\mbox{gr}_{\underline{\nu}} \mathcal{O}_S$ is graded isomorphic
to $\C[\Gamma_\p]^{\underline{v}}$ (see (\ref{graded})). It
follows that:
\[
P^{\underline{\nu}} _{ (S,0)} = P^{\underline{\nu}} _{
(Z^{\Gamma_\p},0)}
\]
for $\Gamma_\p$ the semigroup associated to $(S,0)$ by the
quasi-ordinary projection $\p$ (see section \ref{qo}).

We deduce from Proposition \ref{naive} and Theorem
\ref{caso-torico} that the  Poincaré series expresses as the
quotient of two polynomials in $\C[\underline{t}]$, which are
products of cyclotomic terms.
\begin{corollary}
We have that
\begin{equation} \label{poincare}
P^{\underline{\nu}} _{ (S,0)} = \Phi^{\underline{v}}
(P_{\Gamma_\p}) =  \prod_{j=1}^d \frac{1}{1-
\underline{t}^{\langle \underline{v} , e_j \rangle }}
\prod_{i=1}^g \frac{1- \underline{t}^{\langle \underline{v}, n_i
\g_i \rangle } }{1 - \underline{t}^{\langle \underline{v},
\gamma_i \rangle}} \in \C(\underline{t}).
\end{equation}
\end{corollary}

We exhibit the Poincaré series as an integral with respect to the
Euler characteristic of the projectivization of the space of the
local ring $\mathcal{O}_S$. First we define the $k$-jet space of
functions in $R= \mathcal{O}_S$, relative to the divisorial
valuation $\nu_1$ by (\ref{val-jets}). Then define cylindrical
sets in $\P R$ and their Euler characteristic by using Definition
\ref{cylindric}. Then we have that:
\begin{theorem} \label{euler-q}
\[
P_{(S, 0) }^{\underline{\nu}} = \int_{\P \mathcal{O}_S}
t^{\underline{\nu}} d \chi.
\]
\end{theorem}
\textbf{Proof:} Notice that $\Phi (\P R) = \{ t^{\underline{a}}
\mid \underline{a} \in \langle \underline{\nu}, \Gamma_\p \rangle
\}$. If $\underline{a} \in \langle \underline{\nu}, \Gamma_\p
\rangle $ then (\ref{asteris}) hold.

By remark \ref{function} for any $\g \in \Gamma_\p$ we fix 
an element  $\f_\g \in R$ such that $\f_\g = X^{\g} \epsilon_\g$,
with $\epsilon_\g$ a unit in the integral closure of $R$. Then the
set $ \{ \f_\g \mid \langle \underline{v}, \g \rangle =
\underline{a} \} $ has exactly $\dim_\C J(\underline{a}) $
elements for any $\underline{a}$, by Proposition \ref{key3}. These
elements have linearly independent images in the  $a_1$-jet space
$\mathcal{J}_R^{a_1}$, which span a $\C$-vector space
$L(\underline{a})$ such that
\[ \Phi^{-1} (t^{\underline{a}}) =\p_{a_1}^{-1} ( \P L
(\underline{a}) ).
\]
Then it follows that $\Phi^{-1} (t^{\underline{a}}) $ is a
cylinder and
\[
\chi ( \Phi^{-1} (t^{\underline{a}})) = \chi ( \P
L(\underline{a})) = \chi ( \P J(\underline{a}) ).
\]
We deduce that $t^{\nu}$ is integrable and:
\[
\int_{\P R} t^{\underline{\nu}} d \chi = \sum_{\underline{a} \in
\langle \underline{\nu}, \Lambda \rangle } \chi ( \P
J(\underline{a}) ) t^{\underline{a}} = P_{(Z^\Gamma_\p, 0)
}^{\underline{\nu}} . \quad \Box
\]

\section{Essential divisors} \label{divi}

In this section we recall the definition and characterization of
essential divisors on germs of toric varieties and  of
quasi-ordinary hypersurface singularities. Using the divisorial
valuations associated with essential divisors we define a Poincaré
series as an analytic invariant of the germ section \ref{pser}.

We give some definitions in the algebraic case. Let $\psi_i : W_i
\rightarrow V$ be two modifications, with $W_i$ normal for
$i=1,2$. If $E \subset W_1$ is a divisor then $\psi_2^{-1} \circ
\psi_1 $ is defined on a dense open subset $E^\circ$ of $E$ and we
say then that the closure of $\psi_2^{-1} \circ \psi_1 ( E^\circ)
$ is the {\em center} of the divisor $E$ on $W_2$. We identify
frequently exceptional divisor with the associated divisorial
valuation.

For a variety $V$, a resolution of singularities $\f: V'
\rightarrow V$ is a modification, with $V'$ smooth and such that
$\f$ is an isomorphism over $V - \mbox{Sing}V$. If $E$ is an
exceptional divisor of $\f$ then $\f(E) \subset \mbox{\rm Sing}
V$, i.e.,  $E$ is an {\em exceptional divisor over the singular
locus} of $V$. If $\f(E)$ is contracted to a point $x \in
\mbox{\rm Sing} V$ we say that $E$ is an {\em exceptional divisor
over} $x$. An exceptional divisor $E$ over the singular locus of
$V$ (resp. exceptional divisor over a point $x \in \mbox{Sing} V$)
is {\em essential} (resp. {\em essential over $x$}) if for all
resolution $\f: V' \rightarrow V$ the center of $E$ on $V'$ is an
irreducible component of $\f^{-1} ( \mbox{\rm Sing} V )$ (resp. an
irreducible component of $\f^{-1} (x)$).

In the normal toric case the essential divisors admit a
combinatorial description shown by Bouvier, for essential divisors
over $0$-dimensional orbits (see \cite{Bouvier}), and by Ishii and
Kollar in general (see \cite{IK}).

Let $Z_{\s}$ be an affine toric variety defined by a rational
strictly convex cone $\s \subset N_\R$. We suppose that $Z_{\s}$
is not smooth. Recall that the singular locus of $Z_{\s}$ is the
union of orbits corresponding to non regular faces $\t $ of $\s$.
We denote by $\preceq_{\s}$ the partial order on $N$ defined by:
\[
v \preceq_{\s} v' \Leftrightarrow v' \in v + \s.
\]
\begin{proposition} (see \cite{Bouvier}, \cite{IK})
\begin{enumerate}
\item The essential divisors of the toric variety $Z_{\s}$ are the
divisors corresponding to the minimal elements with respect to the
order $\preceq_{\s}$ in the set
\[
\mathcal{D}_{\s} :=\bigcup_{\t \; \mbox{\rm \small non
regular}}^{\t< \s} \stackrel{\circ}{\t} \cap N.\]

\item The essential divisors over $ o_\s$ are the divisors
corresponding to the minimal elements with respect to the order
$\preceq_{\s}$ in the set $\stackrel{\circ}{\s} \cap N$.
\end{enumerate}
\end{proposition}

By passing it is convenient to give the characterization of
essential divisors of an affine toric variety of the form $(V,0)=
(Z^\Lambda,0)$ for $\Lambda$ a semigroup as in section
\ref{toric}. We denote by $n:\bar{V} \rightarrow V$ the normalization map.
We have that the essential divisors of $V =
Z^\Lambda$ over $0$ coincide with the essential divisors of its
normalization $\bar{V} = Z_{\s}$ over the closed orbit. The key
point to determine the essential divisors over the singular locus
is that the preimage of the singular locus of $V$ is union of
orbits (the normalization map $n$  been equivariant). Some of these
orbits are contained in the singular locus of the normalization
$\bar{V}= Z_{\s, N}$, while other orbits are formed of non
singular points. If $\t$ is a regular face of $\s$ we denote by
$v_{\t}$ the sum of the primitive vectors of $N$ lying in the
edges of $\t$. The following description is due to Ishii (see
\cite{Ishii-crelle}).

\begin{proposition} \label{pretoric}
Let us define
\[
\mathcal{D}'  := \{ v_\t \mid \t \mbox{ regular and } \mbox{orb}(\t)
\subset n^{-1} ( \mbox{\rm Sing} V) \}.
\]
Then the essential divisors of $V$ over its singular locus are the
divisors corresponding to the minimal elements with respect to the
order $\preceq_{\s}$ in the set $ \mathcal{D}_\s \cup
{\mathcal{D}'}$. Moreover, the elements in $\mathcal{D}'$ are  all
minimal in $\mathcal{D}_\s \cup \mathcal{D}'$.
\end{proposition}

The notion of essential divisor extend to the algebroid situation
(see \cite{Ishii-fourier}). We describe precisely the essential
divisors associated to an irreducible germ of quasi-ordinary
hypersurface. Ishii proved that the essential divisors of $(S,0)$
over the origin coincide with the essential divisors of the
normalization $(\bar{S},0) $
over the origin (see
\cite{Ishii-fourier}), where $(\bar{S},0) = (Z_{\s, N}, o_\s)$ by
Proposition \ref{normalization}. For the essential divisors over
the singular locus, the first author proved the following
description, which is based on  previous results of  Ishii in
\cite{Ishii-fourier} and \cite{Ishii-crelle}. The key point as in
the affine toric case is that the preimage by the normalization
map of the singular locus of $S$ is the germ of certain union of
orbits closures.

\begin{proposition} \label{essential}
(See  \cite{GP-nash}) The statement of Proposition \ref{pretoric}
is true when we replace the toric variety $V$ by the germ $(S,0)$.
\end{proposition}

\subsection{The essential matrix} \label{matrix}

Using the characterization of the essential valuations of $(S,0)$
we extract some more combinatorial information on the
characteristic data associated with the quasi-ordinary polynomial
$f$. The content of this section is a technical part needed to
prove the main results of the paper in the next section.

 We consider the
ring extension $\C[[\underline{x}]] \rightarrow \C[[\s^\vee \cap
M]]$, which maps $x_i \mapsto X^{e_i}$, and was introduced in
Section \ref{qo}.
\begin{remark} \label{convention}
We can identify $N_0$ with $\Z^{d}$ by taking coordinates with
respect to the dual basis of $\{ e_i \}_{i=1}^d$. By this
identification the vectors in  $\s \cap N $ as those vectors $v$
in $\Z^d_{\geq 0}$ such that $\langle v, \g_j \rangle \in \Z$ for
$j=1, \dots, g$. We denote by $u_i \in \s$ the primitive vector of
the lattice $N$ which is orthogonal to $e_1, \dots, \hat{e}_i,
\dots, e_d$ (where the term $\hat{e}_i$ is removed from the list).
Notice that the vector $u_i$ is an integral multiple of the
$i^{th}$ vector of the canonical basis of $\Z^d$, for $i=1, \dots,
d$. In particular if $c< d$ then $u_{c+1}, \dots, u_{d}$ are the
vectors in the canonical basis of $\Z^d$.
\end{remark}

\begin{lemma}  \label{sing}
The preimage of the singular locus by the normalization map is
formed by  orbit closures of codimension one or two.

The orbit closures of codimension one correspond to edges of $\s$
of the form $u_i \R_{\geq 0}$ such that $Z_i$ is an irreducible
component of the singular locus of $(S,0)$. The essential divisors
over the codimension one components of $S$ correspond to the
lattice vectors $u_1, \dots, u_{s_1}$.

The orbit closures of codimension two correspond to two
dimensional faces of $\s$ of the form $\t_{i,j} := \mbox{pos}(u_i,
u_j) $ such that $Z_{i,j}$ is an irreducible component of the
singular locus of $(S,0)$. There are $n_{g} -1$ essential divisors
over each codimension two component of $\mbox{\rm Sing} S$, where
$n_g$ is the last characteristic integer of Definition
\ref{lattices}.
\begin{itemize}
\item  If $n_g
>2$ any two vectors  in $\stackrel{\circ}{\t_{i,j}} \cap N$
defining essential valuations are linearly independent.

\item  If $n_g =2 $ there is exactly one essential valuation over
any irreducible component of the form $Z_{i,j}$, which is defined
by the vector $(u_i + u_j)/2$.
\end{itemize}
\end{lemma}
\textbf{Proof:} We use the notations and results on the structure
of the singular locus, Theorem \ref{singular} (see \cite{Lipman88}
Theorem 7.3 and Corollary 7.3.3) combined with the description of
the essential divisors of Proposition \ref{essential}.

We have a natural isomorphism of pairs
\begin{equation} \label{ng}
(\t_{i,j}, N \cap \t_{i,j}) \rightarrow  (\R^2_{\geq 0}, \{ (v_1,
v_2) \in \Z^2_{\geq 0} / v_1 + v_2 =0 \mod  n_g \} ).
\end{equation}
which is deduced from Remark \ref{convention} by the following
argument:  if $v = (v_1, \dots, v_d) \in \t_{i,j} \cap N$, for
$\t_{i, j}$ a face of $\s$ corresponding to a two dimensional
component, then we have that $\langle v, \g_k \rangle =0$ for
$k=1, \dots, g-1$ and   $\langle v, \g_g \rangle = \frac{1}{n_g}
(v_i + v_j)$, since by Theorem \ref{singular} and \ref{rel-semi}
the coordinates $i$ and $j$ of $\g_g$ are equal to $1/n_g$.

It follows immediately from this and Proposition \ref{essential}
that there are exactly $n_g -1$ minimal vectors in the set
$\stackrel{\circ}{\t}_{i,j} \cap N$ with respect to the order
$\preceq_{\s}$, and if $n_g >2$ any two of these vectors are
linearly independent.
 If $n_g =2$ then there is only one minimal
vector $v$ such that  $v_i = v_j =1$ and $v_k =0$ for $k\ne i, j$
and then it follows that $v = (1/2) (u_i + u_j)$. Notice that in
the quasi-ordinary case the cone $\t_{i,j}$, corresponding to a
codimension two component of the singular locus of $S$, is never
regular for the lattice $N$. $\Box$

\begin{definition} \label{esse}
We denote by 
$w_1, \dots, w_{s_1}$  (resp.  $w_{s_1 +1}, \dots, w_{s_1 + s_2}$)   
the sequence of integral vectors in $\s \cap N$ corresponding to
essential valuations of the germ $(S,0)$ centered
on codimension one (resp. two) components of $\mbox{Sing} S$. 
We denote by $w_{s_1 +s_2 +1 }, \dots, w_{ s_0 + s_1 + s_2}$
 the sequence of integral vectors in $\stackrel{\circ}{\s} \cap N$ corresponding to
essential valuations of $(S,0)$ over $0$. We denote by $p$ the number $p= s_0 + s_1 + s_2$. 
We denote by  $\underline{w}$ the following sequence of essencial valuations:
\begin{equation}\label{essentialval}
\underline{w} := (w_1, \dots, w_{p}).
\end{equation}
\end{definition}
In the case $d=2$ notice that if
the singularity is isolated, $s_1 =0$, and  the notions of
essential valuations over the singular locus and over the origin
coincide.

By Theorem \ref{singular} and Lemma \ref{sing} we have that:
\begin{remark} \label{s4}
Notice that $s_0 \geq 1$ by definition of essential valuations of $(S,0)$ over $0$. 
The number $s_1$ equals to the number of codimension one
components of the $\mbox{Sing}\, S$, while the number $s_2$ is
greater or equal to the number $ \left(
\begin{array}{c}
 c - s_1
\\
2
\end{array}
\right) $
 of codimension two irreducible
components of $\mbox{Sing} S$.

We have that $s_2 =1$ if and only if $c- s_1 = n_g = 2$. If $s_2
=1$ then the following equality holds
\[
 \langle w_{s_1 + 1}, e_{c-1} \rangle =  \langle w_{s_1 + 1}, e_{c}
 \rangle= \langle w_{s_1 + 1}, \l_g \rangle=
 1.\]
\end{remark}

 We denote by $\langle \underline{w}, e_j \rangle $ the vector $(
\langle w_1, e_j \rangle, \dots,  \langle w_p, e_j \rangle)$.
\begin{definition}
The matrix
\[\mathcal{M}_f^{\underline{w}} := ( \langle w_i,
e_j \rangle)_{i=1, \dots, p} ^{j= 1, \dots, d}  \] has non
negative integral coefficients and is called the essential matrix
associated with the quasi-ordinary polynomial $f$.
\end{definition}

\begin{remark} \label{linear}
The linear map $\psi_{\underline{w}}: M_{\Q} \rightarrow \Q^{p}$,
which maps
\[
e_j \mapsto \langle \underline{w}, e_j \rangle := (\langle w_1,
e_j \rangle, \dots,  \langle w_p, e_j \rangle), \mbox{ for } j=1,
\dots,d.
\]
has associated matrix $\mathcal{M}_f^{\underline{w}}$ with respect
to the canonical basis.
\end{remark}

The following  technical Lemma implies that the matrix
$\mathcal{M}_f^{\underline{w}}$ has the following block structure:
\begin{center}
\begin{tabular}{|c|c|c|c|}
  \hline
  & ${s_1} $&  ${c- s_1}$ & ${d- c} $
  \\
  \hline
  $s_1 $ & $D$ & $\underline{0}$ & $\underline{0}$ \\
  $s_2 $ &$ \underline{0} $&  $B$& $ \underline{0} $\\
  $s_0 $  & $\underline{*} $ &  $\underline{*} $ &  $\underline{1}$\\
\hline
\end{tabular}
\end{center}
where the first line  and column indicates the size of the blocks.
The submatrix $D$ is diagonal of maximal rank $s_1$, $B$ is
generically of maximal rank $c-s_1$, $\underline{*}$ means non
zero entry and $\underline{1}$ or $\underline{0}$ means all
entries are one or zero, respectively. More precisely we have
that:
\begin{proposition} \label{suficiente}
If $d> 2$ the following statements hold.
\begin{enumerate}
\item \label{s1} The submatrix $ D:= ( \langle w_i, e_j \rangle
)_{1\leq i, j \leq s_1}$ is diagonal with non zero determinant.

\item \label{s2} The submatrices $(\langle w_i, e_j \rangle) _{s_1
+1 \leq i \leq s_1 + s_2}^{ 1\leq j \leq s_1}$, $(\langle w_i, e_j
\rangle) _{s_1 +1 \leq i \leq s_1 + s_2}^{c< j\leq d} $ and
$(\langle w_i, e_j \rangle) _{1 \leq i \leq s_1}^{ s_1 < j \leq d}
$ vanish.

\item \label{s3} The submatrix $(\langle w_i, e_j \rangle) _ {s_1+
s_2 +1\leq i \leq p}^{ c+1 \leq j \leq d}$ have all entries equal
to one.

\item \label{s5} If $c -s_1 \geq 2$ then the submatrix $B:= (
\langle w_i, e_j \rangle )_{s_1 + 1\leq i \leq s_1 + s_2}^{ s_1 +
1 \leq j \leq c} $ is of rank $c- s_1$ if $s_2 \ne 1$.

\item \label{s6} The submatrix $( \langle w_i, e_j \rangle
)_{1\leq i \leq p}^{ 1\leq j \leq c} $ is of rank $c$ if $s_2 \ne
1$.
\end{enumerate}
\end{proposition}
\textbf{Proof:} By Lemma \ref{sing} we have that $w_1 = u_1,\dots,
w_ {s_1} = u_{s_1} $. Then the first assertion follows from this
since  $ \langle u_i, e_j \rangle = 0 \Leftrightarrow i \ne j,
\mbox{  for } 1 \leq i, j \leq d$.

By Lemma \ref{sing} we have that if $s_1 +1 \leq  i \leq  s_1 +
s_2$ then there are $s_1 < j_1 < j_2 < c$ such that $w_i \in
\mbox{\rm pos} (u_{j_1}, u_{j_2})$. Then assertion {\em \ref{s2}}
follows from:
\begin{equation} \label{caso}
\langle w_i , e_j \rangle =0 \Leftrightarrow j \ne j_1, j_2,
\mbox{ for } 1 \leq j \leq d.
\end{equation}

Notice that if $c< d$ the vectors $u_{c+1}, \dots, u_d$ are part
of the canonical vectors in $\Z^d \supset N$ by Remark
\ref{convention}. The semigroup $\s \cap N$ splits in the form
$(\s' \cap N')\times \bigoplus_{k=c+1}^d \Z u_k$, where $\s'$
(resp. $N'$) is the intersection of the cone $\s$ (resp. the
lattice $N$) with the subspace spanned by $\{ u_1, \dots, u_c \}$.
This implies that the vectors $w$ corresponding to essential
divisors in $\s$ are of the form
 $w = w' + u_{c+1} + \cdots + u_d$ where $w'$ corresponds
 to a vector defining an essential divisor in $\s'$.
 This shows statement {\em \ref{s3}}.

For assertion {\em \ref{s5}} we have that  if $s_2 \geq 2 = c-
s_1$ the vectors $w_{s_1 +1}, w_{s_1+2}$, viewed as forms on $M$,
have linearly independent restrictions to the sublattice spanned
by $e_{c-1}, e_c$.

Otherwise if $c- s_1 > 2$ then the number $s_2$ is greater equal
than the number of codimension two components of $\mbox{\rm Sing}
S$, which is $\geq c- s_1$.  It is easy to check the result using
(\ref{caso}).

For statement {\em \ref{s6}} it is enough to verify the result in
the case $c- s_1=1$ (and then $s_2=0$) by {\em \ref{s1}} and {\em
\ref{s5}}. The forms $w_1, \dots, w_{c-1}, w_{c}$, restricted to
the subspace generated by $e_1, \dots, e_c$, are linearly
independent by assertion {\em \ref{s1}}, since $w_c \in
\stackrel{\circ}{\s}$ takes non zero value on $e_c$, while $w_1,
\dots, w_{c-1}$ vanish on $e_{c}$. $\Box$

\section{The Poincaré series as an analytical invariant of the singularity}
\label{pser}

This definition of Poincaré series of the toric or quasi-ordinary
hypersurface singularity depends obviously on the set of vectors
$v_1, \dots, v_r$. However we get rid of this dependency  by
considering the Poincaré series associated with the essential
valuations, which in this case are monomial valuations.

In the toric case we have:
\begin{definition}
The Poincaré series $P_{(Z^\Lambda, 0)}$ of the germ $(Z^\Lambda,
0)$ is the Poincaré series associated to the sequence of essential
valuations  of $(Z^{\Lambda}, 0)$ over the singular locus together
with the essential valuations of $(Z^{\Lambda}, 0)$ over the
origin.
\end{definition}

This  series is well defined since any essential valuation over
the origin corresponds to a primitive vector in
$\stackrel{\circ}{\s} \cap N$, hence (\ref{chi}) holds. Notice
that $P_{(Z^\Lambda, 0)} \in \Z[[ \underline{t}]]$ is an analytic
invariant of the germ $(Z^\Lambda, 0)$  up to permutation of the
variables $t_1, \dots, t_r$.

If $(S,0)$ is an analytically irreducible germ of quasi-ordinary
hypersurface singularity we define an intrinsic notion of Poincaré
series:
\begin{definition} \label{pserie}
The  Poincaré series   $P_{(S, 0)}$   of the germ $(S, 0)$ is the
Poincaré series associated to the sequence (\ref{essentialval}) of essential valuations
of $(S, 0)$ over the singular locus and the  essential valuations
of $(S, 0)$ over the origin.
\end{definition}

\begin{remark} \label{esse2}
The indeterminates $\underline{t} =( t_1, \dots, t_p)$ in
Definition \ref{pserie} are distinguished in groups as follows: If
$d>2$ then $t_1, \dots, t_{s_1}$ (resp. $t_{s_1 +1}, \dots,
t_{s_2}$) corresponding to the essential valuations
$\underline{w}= (w_1, \dots, w_p)$ over codimension  one (resp.
two) components of $\mbox{Sing} S$, while $t_{s_1+ s_2 +1}, \dots,
t_p$ corresponding to those essential valuations over the origin.
If $d=2$ we consider only two groups as in the case of essential
valuations (see the remarks after Definition \ref{esse}). The
Poincaré series of  $(S,0)$ is an analytic invariant of the germ
$(S,0)$, up to permutations of $t_1, \dots, t_p$ respecting these
groups of indeterminates introduced before.
\end{remark}

\begin{example} \label{quadratic}
We consider the quadratic cones $y^2 -x_1 x_2 =0$  in $\C^{d+1}$,
for $d \geq 2$.
\end{example}

\begin{enumerate}
\item If $d=2$ then there is only one
 essential valuation over the
singular locus, which is reduced to the origin. It follows that
the associated Poincaré series depends only on one variable $t$
and we find that:
\[
P_{(S,0)} = (1-t^2)(1-t)^{-3}.
\]

 \item If $d> 2$ then there are two essential valuations, one
over the singular locus, and the other over the origin,
corresponding respectively to the variables $t_1, t_2$ in the
associated Poincaré series
\[
P_{(S,0)} = (1-t_1^2 t_2^2)(1- t_1 t_2)^{-3} (1-t_2)^{-d+ 2}.
\]

\end{enumerate}

The main result of this paper is the following:
\begin{theorem} \label{Pseries}
The Poincaré series of an irreducible germ $(S, 0)$ of
quasi-ordinary hypersurface determines (and it is determined by) the
sequence of normalized characteristic monomials of $(S,0)$.
\end{theorem}

We analyze the terms appearing in the Poincaré series $P_{(S,0)}$
when the sum is obtained as in (\ref{poincare}), from a normalized
quasi-ordinary projection and the indeterminates $\underline{t} =
(t_1, \dots, t_p)$ are labelled as in Remark \ref{esse2}.

\begin{lemma} \label{cyclotomic} We suppose $d > 2$.
If the expression   (\ref{poincare}) is deduced from a normalized
quasi-ordinary projection of a quasi-ordinary hypersurface of
dimension $>2$  then there is no cancellation between cyclotomic
terms in the numerator and the denominator. Moreover, if we denote
by $c$ the equisingular dimension of $(S,0)$ then we have that:
\begin{enumerate}
\item  \label{c1} The dimension $d$ of $(S,0)$ is the number of
cyclotomic factors appearing in the denominator minus the number
of cyclotomic factors appearing in the numerator in the expression
(\ref{poincare}), counted with multiplicities.

\item  \label{c2} The number $g$  of normalized characteristic
exponents  is the number of cyclotomic factors in the numerator.

 \item \label{c3} If $s_2= 0$ the term $(1 - t_{s_1 + s_2 +1} \cdots t_{p})$ appears in the denominator
of (\ref{poincare}) with multiplicity equal to
\[
\left\{
\begin{array}{lclcl}
d-c & \mbox{ if } & c = d  & \mbox{ or } & \langle \underline{w},
e_c \rangle \ne \langle \underline{w}, e_d \rangle,
\\
d - c +1 & \mbox{ if } & c \ne d & \mbox{ and } &  \langle
\underline{w}, e_c \rangle = \langle \underline{w}, e_d \rangle.
\end{array}
\right.
\]

\item \label{c4} If $s_2 =1 $ the term $(1 - t_{s_1 + s_2 +1}
\cdots t_{p})$ appears with multiplicity $d- c$ and the term $(1 -
t_{s_1 +1} \cdots t_p)$ appears with multiplicity equal to
\[
\left\{
\begin{array}{lcl}
3 & \mbox{ if }  & g = 1 \mbox{ and } \g_1 = (1/2, 1/2, 0,
\dots,0),
\\
2 & & \mbox{ otherwise. }
\end{array}
\right.
\]

\item \label{c5} If $s_2 \geq 2$ then the term $(1 - t_{s_1 + s_2
+1} \cdots t_{p})$ appears with multiplicity $d- c$.
\end{enumerate}
In all the cases the cyclotomic terms not mentioned above appear
without multiplicities.
\end{lemma}
\textbf{Proof:} The assertion one is trivial. By (\ref{preforma})
we have that
\[
\g_1 < n_1 \g_1 < \g_2 < n_2 \g_2 < \cdots < \g_g < n_g \g_g,
\]
(where $<$ above means $\ne$ and $\leq$ coordinate-wise).  Since
$w_{s_1+s_2+1} \in \stackrel{\circ}{\s}$ we deduce that
\begin{equation} \label{desi}
\langle \underline{w}, \g_1 \rangle < \langle \underline{w}, n_1
\g_1\rangle <  \cdots < \langle \underline{w},\g_g \rangle<
\langle \underline{w}, n_g \g_g\rangle.
\end{equation}
 A ray of the form $\R_{\geq 0} \g_j$ may contain  only the
vector $e_1$ in $\{e_1, \dots, e_d\}$, and if this is the case,
since the quasi-ordinary projection is normalized, we have that
$e_1 < \g_1 = \l_1$ and then
\begin{equation}\label{desi2}
\langle \underline{w}, e_1 \rangle < \langle \underline{w}, \g_1
\rangle.
\end{equation}
A cancellation  between terms in the numerator and denominator of
(\ref{poincare}) implies the existence $1 \leq j \leq g$ and $1
\leq i \leq d$ such that:
\begin{equation} \label{am}
\langle \underline{w}, n_j \g_j \rangle = \langle \underline{w},
e_i \rangle.
\end{equation}
By definition $\g_j $ is an element in the $\Q$-vector space
generated by $e_1, \dots, e_c$. Then Proposition \ref{suficiente}
number {\em 3}, implies that any such $i$ has to be $\leq c$,
since some positive power of $t_k$, for $1 \leq k \leq s_1 +s_2$
divides $\underline{t}^{\langle \underline{w}, \g_j \rangle}$. We
distinguish two cases:

- If $s_2 \ne 1$ there is no possible cancellation between terms
in the numerator and denominator of (\ref{poincare})
 since the matrix
$(\langle w_i, e_j \rangle)_{1\leq i \leq p}^{1 \leq j \leq c}$
 is of maximal rank $c$ by Proposition \ref{suficiente}.

- If  $s_2 =1 $ and if $1 \leq j < g$ and $1 \leq i < c -1$
equality in (\ref{am}) may not occur: the previous discussion
implies that $e_i = e_1$ but then (\ref{desi}) and (\ref{desi2})
give a contradiction. By remark \ref{s4} number {\em 2} we have
that $n_g =2$ and
\[
\underline{t}^{\langle \underline{w}, e_{c-1} \rangle} =
\underline{t}^{\langle \underline{w}, e_{c} \rangle} = t_{s_1+1}
\cdots t_{p}.
\]
It is clear that this term is not equal to $\underline{t}^{\langle
\underline{w}, n_j \g_j \rangle}$ for $j=1, \dots, g-1$ while the
exponent of $t_{s_1 +1}$ in $\underline{t}^{\langle \underline{w},
2 \g_g \rangle}$ is equal to two by remark \ref{s4} number {\em
2}. No cancellation is possible in this case. Example
\ref{quadratic} is a particular case of this situation.

Now assertion {\em \ref{c2}} is clear. Statement {\em \ref{c3}}
follows from Proposition \ref{suficiente}, number {\em \ref{s3}},
taking into account that if $d> c$ then the column vectors
$\langle \underline{w}, e_c \rangle $, $\langle \underline{w}, e_d
\rangle$ of ${\mathcal M}_f^{\underline{w}}$ may be equal. The
statement {\em \ref{c5}} is similar. The first statement in
assertion {\em \ref{c4}} follows from Proposition
\ref{suficiente}, number {\em \ref{s3}} and {\em \ref{s5}}. The
second statement follows from Remark \ref{s4} and a simple
computation.

Finally, no other cyclotomic term appears with multiplicity $>1$
since Proposition \ref{suficiente} number {\em \ref{s6}} implies
that the column vectors $\langle \underline{w}, e_j \rangle _{j=1,
\dots, c}$ are linearly independent and determine $\langle
\underline{w}, \l_i \rangle$  by linearity over $\Q$ (see Remark
\ref{linear}). $\Box$

\textbf{Proof of Theorem \ref{Pseries}} After eventual
simplification we have that by Lemma \ref{cyclotomic} the Poincaré
series of $S$ expresses uniquely in the form:
\[
P_{(S,0)} = \frac{ (1 - \underline{t}^{\b_1}) \cdots  (1 -
\underline{t}^{\b_{g}}) }{ (1 - \underline{t}^{\a_1}) \cdots  (1 -
\underline{t}^{\a_{d+g}})
 }
\]
where $d$ is the dimension of $S$ and $g$ is the number of
normalized characteristic exponents. By (\ref{desi}) we can
relabel if necessary in order to have $\b_1 < \dots < \b_g$. It
may happen that the ray $\b_1 \R_{\geq 0}$ contains several of the
$\a_i$ or at most one. We denote in the first case by $\a_1$ the
smallest (coordinate-wise) and by $\a_{d+1}$ the following one
with this property. In the second case we denote by $\a_{d+1}$ the
unique exponent with this property. By Lemma \ref{cyclotomic} we
have that $n_1 \a_{d+1} = \b_1$ for $n_1 >1$. We repeat the
procedure relabelling if necessary in order to have that $n_i
\a_{d+i} = \b_i$ for $i=1, \dots, g$. Notice that the integers
$n_1, \dots, n_g$ are characteristic integers associated with the
normalized characteristic exponents of $(S,0)$.

We deal first with the case $d>2$. The Poincaré series distinguish
the $s_2$ indeterminates corresponding to the essential divisors
over a  codimension two component of the singular locus of $(S,0)$
(see Remark \ref{esse2}).

- If $s_2 \geq 2$ we have that by Lemma \ref{cyclotomic} the term
$1 - t_{s_1 + s_2 +1} \cdots t_p$ appears in the denominator of
$P_{(S,0)}$ with multiplicity $d-c$. This determines $c$. Notice
that by Remark \ref{s4} the case $c -s_1 = n_g =2$ is excluded. We
relabel then the exponents $\a_1, \dots, \a_d$ in such a way that
\begin{equation} \label{index}
\a_{d-c+1} = \cdots = \a_d = (\stackrel{c}{\overbrace{0, \dots,
0}},1, \dots, 1).\end{equation} The matrix $\mathcal{M}$, which
has  columns  the vectors  $\a_1, \dots, \a_c$, is equal to
$(\langle w_i, e_j \rangle)_{1 \leq i \leq p}^{1 \leq j \leq c}$
up to a permutation of the columns, by Lemma \ref{cyclotomic}. By
Proposition \ref{suficiente} number {\em \ref{s6}} the matrix
$\mathcal{M}$ is of maximal rank $c$. The linear equation
$\mathcal{M} \l = \a_{d +j} $ has exactly one solution, which give
the first $c$ coordinates of $\g_j$, eventually under a fixed
permutation, for $j=1, \dots, g$. We recover the vectors $\g_1,
\dots, \g_g$ in the minimal sequence of generators of the
semigroup $\Gamma_\p \subset \Q^d$ by adding to each solution
$d-c$ zeroes in order to have a $d$-tuple and relabelling the
coordinates  of the resulting vectors in such a way that
(\ref{lex}) holds. Formula (\ref{rel-semi}) determines from these
data the normalized sequence of characteristic exponents.

- If $s_2 =0$ then  by Lemma \ref{cyclotomic} the term $1 - t_{s_1
+1} \cdots t_p$ appears in the denominator of $P_{(S,0)}$ with
multiplicity $d- c'$ where $c' \in \{ c-1, c \}$. We relabel the
exponents $\a_1, \dots, \a_d$ in such a way that those
corresponding to $t_{s_1 +1} \cdots t_p$ appear with the highest
index (as in (\ref{index}) above). By Lemma \ref{cyclotomic} the
matrix $\mathcal{M}$, with columns the vectors $\a_1, \dots,
\a_{c'}$ is equal to $(\langle w_i, e_j \rangle)_{1 \leq i \leq
p}^{1 \leq j \leq c'}$, up to a permutation of the columns. By
Proposition \ref{suficiente} number {\em \ref{s6}} the matrix
$\mathcal{M}$ is of maximal rank $c'$.

In order to determine $c$ we study the linear equation
$\mathcal{M} \l= \a_{d +g} $. If this equation has no solution
then $c= c' +1$ and then we replace $\mathcal{M}$ by the matrix,
denoted with the same letter, with columns the vectors $\a_1,
\dots, \a_{c'+1}$. In this case this matrix is equal to $(\langle
w_i, e_j \rangle)_{1 \leq i \leq p}^{1 \leq j \leq c'+1}$, up to a
permutation of the columns, and it is of rank $c'+1 =c$, by
Proposition \ref{suficiente} number {\em \ref{s6}}. Then in both
cases the linear system $\mathcal{M} \l = \a_{d+j} $ has exactly
one solution
 for $j=1, \dots,g$ from which we recover the
normalized characteristic exponents as in the previous case.

- If $s_2 =1$ then   the term $1 - t_{s_1 + 2} \cdots t_p$ appears
in the denominator of $P_{(S,0)}$ with multiplicity $d-c$ and the
term $(1 - t_{s_1 +1} \cdots t_p)$ appears also with multiplicity
$2$ or $3$.
 We relabel
the exponents $\a_1, \dots, \a_d$ in such a way that those
corresponding to the monomial $t_{s_1 +2} \cdots t_p$ appear with
the highest index (as in (\ref{index}) above) followed by those
exponents corresponding to the monomial $t_{s_1 +1} \cdots t_p$.
Notice that we can distinguish these monomials in the Poincaré
series by remark \ref{esse2}. We denote by $\xi'$ the $s_1$-tuple
obtained from $\xi \in \Q^p$ by keeping only the first $s_1$
coordinates. By Lemma \ref{cyclotomic} the matrix $\mathcal{M}$,
with columns the vectors $\a_1', \dots, \a_{s_1}'$ is equal to
$(\langle w_i, e_j \rangle)_{1 \leq i \leq s_1}^{1 \leq j \leq
s_1}$, up to a permutation of the columns. By Proposition
\ref{suficiente} number {\em \ref{s1}} we have that $\mathcal{M}$
is then of maximal rank $s_1$. The linear equation $\mathcal{M} \l
= \a_{d +j}' $ has exactly one solution, for $j=1, \dots, g-1$,
which give the first $c-2= s_1$ coordinates of $\g_j'$, eventually
under a fixed permutation, for $j=1, \dots, g$. We form $d$-tuples
from these terms putting zeroes in the new coordinates for $j=1,
\dots, g-1$. For $j=g$ we put two new coordinates equal to $1/2$
and the remaining ones equal to zero (we use Remark \ref{s4},
Theorem \ref{singular} and Formula \ref{rel-semi}). Then the proof
finish as in case $s_2 \geq 2$.

Finally we deal with the surface case. If the singularity is
normal it is isolated hence it follows that $s_1 =0$,  $c= 2$ and
the essential valuations over the singular locus coincide with the
essential valuations over the origin. In this case  the
singularity is isomorphic to $y^n- x_1x_2$ by \cite{Lipman88}. It
is a $A_{n-1}$ singularity possessing  $n-1$ exceptional divisors
in the minimal resolution, all of them essential divisors
corresponding to integral lattice vectors pair-wise linearly
independent if $n> 2$ (see \cite{Fulton}). If $n=2$ then there is
only one of such exceptional divisor in the minimal resolution and
hence only one essential valuation (see Example \ref{quadratic}).
If the singularity is not normal then its singular locus of $S$
has one or two irreducible components. If we have two components
then the assertion is consequence of Lemma \ref{sing}: the
corresponding vectors $w_1= u_1, w_2 = u_2$ are obviously linearly
independent. Otherwise if there is only one component,
corresponding the vector $w_1= u_1$, then we have that the vector
$w_p \in \stackrel{\circ}{\s}$ is linearly independent with $w_1$.
In any of these cases,with the exception of the quadratic cone
$y^2 = x_1 x_2$, the matrix $\mathcal{M}$ with columns the vectors
$\a_1, \dots, \a_{2}$ is equal to a permutation of the columns of
 $(\langle w_i, e_j \rangle)_{1 \leq i \leq p}^{1 \leq
j \leq 2}$. Since this matrix is non singular we recover the
vectors $\g_1, \dots, \g_g$, up to a fixed permutation of the
coordinates by resolving the linear systems ${\mathcal M} \g =
\a_{2+i}$ for $i=1, \dots, g$. Then we conclude as in the previous
cases. $\Box$.

Notice that the procedure is completely algorithmic once the
reduced form of the Poincaré series is given, and uses only
elementary linear algebra operations.

Let $(S,0) \subset (\C^{d+1}, 0) $ be an irreducible germ of
quasi-ordinary hypersurface singularity defined by a normalized
quasi-ordinary polynomial $f \in \C[[x_1, \dots, x_{d}]][y]$. We
denote by $(S',0)$ the $d'$ dimensional germ obtained from $S$ by
intersecting with $x_{d-d'+1} = \cdots = x_d =0$. The following
corollary is then immediate:

\begin{corollary}\label{equi}
The germ $(S,0)$ is an equisingular deformation of a
quasi-ordinary hypersurface singularity $(S',0)$ of dimension $1
\leq d'\leq d$ in the sense of section \ref{equisingular} if and only if
\[
P_{(S,0)} (t_1, \dots, t_{p}) = \left(\frac{1}{1-t_{s_1+ s_2 +1}
\cdots t_p}\right)^{d-d'} P_{(S',0)} (t_1, \dots, t_{p}),
\]
where the indeterminates $t_1, \dots, t_p$ are distributed in
groups (see Remark \ref{esse2}) in the same form for $(S,0)$ and
$(S',0)$.
\end{corollary}

Notice that the minimal possible  value of $d'$ in Corollary
\ref{equi} is equal to $c$, the dimensional type of $(S,0)$.

\subsection{The Poincaré series and the zeta function} \label{zeta} 

In this section we compare the zeta-function of geometric
monodromy of an irreducible quasi-ordinary polynomial function $f:
(\C^{d+1}, 0) \rightarrow (\C, 0)$ with the Poincaré series of the
germ $f=0$ at $0$, for $d\geq 2$. This zeta function is studied by
McEwan and Némethi in \cite{Nemethi}. See also \cite{NemethiII}
for a description when $f$ is not irreducible.

We suppose that $f \in \C \{ x_1, \dots, x_d \} [y]$ is  a
normalized irreducible quasi-ordinary polynomial, with
characteristic exponents $\l_j = (\l_{j,1}, \dots, \l_{j,d}) \in
\Q^d$ for $j=1, \dots,g$.
\begin{theorem}(Theorem A. McEwan, Némethi\cite{Nemethi}) \\
The zeta function of geometric monodromy of $f: (\C^{d+1}, 0)
\rightarrow (\C, 0)$ is equal to:
$$
\zeta(f)(t)=\zeta(f\mid_{x_{2}=\cdots=x_{d}=0})(t).
$$
More precisely, there are two cases:
\begin{description}
\item[A]If $\lambda_{1,2}\neq 0$ then  we have that $
\zeta(f)(t)=1-t^n$ with $n=deg_{y}(f)$.

\item[B] If $\lambda_{1,2}=0$ and if $i_0$ denotes the maximum
among those index $1 \leq i \leq g$ such that $\l_{i, 2} =0$ let
$h (x_1, y) \in \C\{ x_1, y\} $ be the irreducible series with the
Puiseux expansion $ y = \sum_{\l_i \leq \l_{i_0}}
{x_1}^{\l_{i,1}}$ then we have that $ \zeta(f)(t) =\zeta(h)(t^{
{n}/{n_1 \dots n_{i_0}}})$.
\end{description}
\end{theorem}

Denote by $P $ the result of substituting the indeterminates $t_1,
\dots, t_p$ by $t$ in the Poincaré series $P_{(S,0)}$. Remark that
$P$ coincides with the Poincaré series $P_{(S,0)}^w$ where $w:=
w_1 + \cdots + w_p$, is the sum of all vectors in $\s \cap N$
corresponding to essential valuations. Remark that $w \in
\stackrel{\circ}{\s} \cap N$. Let us write $w=(b_1,\ldots,b_d) \in
\Z^d$, for $b_i=\sum_{j=1}^r w_{i,j}$,  where the coordinates
$w_i= (w_{i,1}, \dots, w_{i,d})$ are defined by Remark
\ref{convention}. Then $P$ expresses as a quotient of products of
cyclotomic polynomials in one variable as follows:
$$
P =
\left(\frac{\prod_{j=1}^{i_0}(1-t^{b_1n_j\gamma_{j,1}})}{\prod_{j=1}^{i_0}(1-t^{b_1\gamma_{j,1}})(1-t^{b_1})}\right)
\prod_{j=i_0+1}^g\left(\frac{1-t^{\langle
w,n_j\gamma_j\rangle}}{1-t^{\langle w,\gamma_j\rangle}}\right)
\prod_{j=2}^d\left(\frac{1}{1-t^{b_j}}\right)
$$

\begin{corollary} \label{fzeta}
In case {\bf A} the function ${\zeta(f)(t^{\frac{b_1}{n}})}$
appears as a cyclotomic factor  in the denominator of $P$.

In case {\bf B} the function $\zeta(h)(t^{\frac{b_1}{deg(h)}})$
appears as products of cyclotomic factors in the cyclotomic
expansion of $P$. More precisely we have that
\end{corollary}

\begin{description}
\item[A] $\lambda_{1,2}\neq 0$, \quad  $ P
=\frac{1}{\zeta(f)(t^{\frac{b_1}{n}})}
\left(\frac{\prod_{j=1}^{i_0}(1-t^{b_1n_j\gamma_{j,1}})}{\prod_{j=1}^{i_0}(1-t^{b_1\gamma_{j,1}})}\right)
\prod_{j=i_0+1}^g\left(\frac{1-t^{\langle
w,n_j\gamma_j\rangle}}{1-t^{\langle w,\gamma_j\rangle}}\right)
\prod_{j=2}^d\left(\frac{1}{1-t^{b_j}}\right) $.
 \item[B]
$\lambda_{1,2}=0$, \quad $ P
=\zeta(h)(t^{\frac{b_1}{deg(h)}})\prod_{j=i_0+1}^g\left(\frac{1-t^{\langle
w,n_j\gamma_j\rangle}}{1-t^{\langle w,\gamma_j\rangle}}\right)
\prod_{j=2}^d\left(\frac{1}{1-t^{b_j}}\right). $
\end{description}

\vspace{.2cm}

Notice that the Poincaré series $P$ encodes the information
provided by the zeta function exactly when $i_0 = g$, i.e., in the
case of an equisingular family of plane branches, i.e., $c=1$.

{\bf Question.} Is it possible to find a notion of Zeta Function
which is more adapted to non isolated singularities and which in
the case of quasi-ordinary hypersurface singularities encodes the
embedded topological type information recovered by the Poincaré
series.

\subsection{Example}
Let us consider a three dimensional germ  $(S,0) \subset
(\mathbb{C}^4,0)$ which in coordinates $(x_1,x_2, x_3, y)$ has
equation $f=0$ where,
\[
f =(y - x_1)^3 - x_1^{11} x_2  (1 + x_3).
\]
The characteristic exponents corresponding to  the quasi-ordinary
projection \[(x_1,x_2, x_3,y) \mapsto (x_1, x_2, x_3)\] are
$\lambda_1=(\frac{1}{3},0,0)$ and
$\lambda_2=(\frac{5}{9},\frac{1}{9},0)$. The generators of the
semigroup $\Gamma_\p$ are $\gamma_1=(\frac{1}{3},0,0)$ and
$\gamma_2=(\frac{11}{9},\frac{1}{9},0)$ and the charateristic
integers $n_1=3$ and $n_2=9$. By Theorem \ref{singular} the
singular locus of $(S,0)$ has only one component, which is of
codimension one, defined by the intersection with $x_1 =0$.

We identify the lattice $N_0$ defined in section \ref{qo} with
$\Z^3$ using Remark \ref{convention}. Then the lattice $\s \cap N$
is a subsemigroup  of $\Z^3_{\geq 0}$ and $N$ is viewed as a
sublattice of $\Z^3$.  We find that the vectors $v_1 = (9,0,0)$,
$v_2 = (3,3,0)$ and $v_3 = (0,0,1)$ are a basis of $N$ (for this
purpose we can use the algorithms introduced by Popescu-Pampu in
\cite{PPP05}). By Lemma \ref{sing} the essential divisor over the
codimension one component of $S$ corresponds to the vector $v_1$,
and the essential divisors over the origin correspond to the
minimal vectors in $\stackrel{\circ}{\s} \cap N$ with respect to
the order $\preceq_{\s}$. We find in this case only one minimal
element $w = (3,3,1)$.

The essential matrix associated with $\bar{w}:= (v, w)$ is:
\[
M_f^{\bar{w}} = \left(%
\begin{array}{ccc}
  9 & 0 & 0 \\
  3 & 3 & 1 \\
\end{array}%
\right)
\]
and the Poincaré series is
\begin{equation} \label{pej}
P_{(S,0)} = \frac{(1- t_{1}^9 t_{2}^3)(1-t_{1}^{99}
t_{2}^{36})}{(1- t_{1}^9 t_{2}^3
)(1-t_{2}^3)(1-t_{2})(1-t_{1}^3t_{2})(1- t_{1}^{11} t_{2}^4)}.
\end{equation}

Notice that the cyclotomic factor ${(1-t_{2}^3t_{1}^9)}$ appears
in both numerator and denominator of (\ref{pej}). This does not
contradicts Lemma \ref{cyclotomic} since the quasi-ordinary
projection is not normalized (see section \ref{qo}). We obtain
then the sort form of the Poincaré series a quotient of products
of cyclotomic polynomials,

\begin{equation} \label{Pej2}
P_{(S,0)} =\frac{1- t_{1}^{99} t_{2}^{36}
}{(1-t_{2})(1-t_{2}^3)(1-t_{1}^3
t_{2}t_{1}^3)(1-t_{1}^{11}t_{2}^4)}.
\end{equation}

We apply now the proof of Theorem \ref{Pseries} as an algorithm
starting from the expression (\ref{Pej2}). We deduce that $d=3$
and $g=1$. Moreover since we have only one relation of the form $(
t_{1}^{11}t_{2}^4)^9 = t_{1}^{99} t_{2}^{36}$ we have that the
first normalized characteristic integer is $9$. By hypothesis we
know that $s_2=0$. The term $1- t_2$ appears in the denominator
with multiplicity one. The matrix $\mathcal{M}$ is equal to
\[
\mathcal{M} = \left(%
\begin{array}{cc}
  3 & 0 \\
  9 & 3 \\
\end{array}%
\right) .
\]
The linear system $\mathcal{M}\l = (11, 4)^t$ has a unique
solution $(\frac{11}{3}, \frac{1}{9})$, therefore the normalized
sequence of characteristic exponents associated to $(S,0)$ in this
case is reduced to $\l_1 ' = (\frac{11}{3}, \frac{1}{9},0)$, the
normalized characteristic sequence associated to $(S,0)$.

\vspace{.3cm}

 {\bf Acknowledgement.} The authors are grateful to
A. Campillo, F. Delgado, A. Lemahieu,  P.Popescu-Pampu and B.
Teissier for helpful suggestions and discussions. The second
author thanks the {\em Universidad Complutense de Madrid} and its
research team {\em Grupo Singular} for their hospitality.

\end{document}